\documentclass[review]{elsarticle}

\usepackage{hyperref}

\journal{Journal of Computational Physics}

\usepackage{color}
\usepackage{graphicx}
\usepackage{ulem}
\usepackage{hyperref}
\usepackage{multirow}
\usepackage{amsmath,amstext,amssymb}
\usepackage{amsthm}
\usepackage{tabularx,longtable,multirow,diagbox}
\usepackage{booktabs}
\newcommand{\jm}[1]{\left[#1\right]}

\newcommand{\R}{\mathbb{R}}

\newcommand{\bn}{\mathbf{n}}

\newcommand{\be}{\beta}

\newcommand{\Ga}{\Gamma}

\newcommand{\Om}{\Omega}
\newcommand{\pa}{\partial}

\newcommand{\eq}[1]{\begin{align}#1\end{align}}

\definecolor{gold}{rgb}{0,0,1}

\newtheorem{remark}{Remark}

\newtheorem{example}{Example}

\numberwithin{equation}{section} \numberwithin{theorem}{section}
\numberwithin{corollary}{section} \numberwithin{lemma}{section}

\graphicspath{{./pictures/}}

\begin{document}
	
	\begin{frontmatter}
		
		%

		\title{XI-DeepONet: An operator learning method for elliptic interface problems}
		
		\author[nju]{Ran Bi} 
		\ead{ranbi@smail.nju.edu.cn}

		\author[ustc]{Jingrun Chen\fnref{myfootnote2}}
		\ead{jingrunchen@ustc.edu.cn}
		\fntext[myfootnote2]{The work of this author was
			partially supported by the NSFC Major Research Plan - Interpretable and General Purpose Next-generation Artificial Intelligence (Nos. 92270001 and 92370205).}
		
		\address[nju]{School of Mathematics,
			Nanjing University, Nanjing 210093, People's Republic of China}
		\address[ustc]{School of Mathematical Sciences and Suzhou Institute for Advanced Research, University of Science and Technology of China, Suzhou 215127, People's Republic of China}
		\author[nju]{Weibing Deng\corref{mycorrespondingauthor}\fnref{myfootnote1}}
		\ead{wbdeng@nju.edu.cn}
		\cortext[mycorrespondingauthor]{Corresponding author}
		\fntext[myfootnote1]{The work of this author was
			partially supported by the NSF of China grant  12171237, and by the
				Ministry of Science and Technology of China grant 2020YFA0713800.}

		\begin{abstract}
		Scientific computing has been an indispensable tool in applied sciences and engineering, where traditional numerical methods are often employed due to their superior accuracy guarantees. However, these methods often encounter challenges when dealing with problems involving complex geometries. Machine learning-based methods, on the other hand, are mesh-free, thus providing a promising alternative. In particular, operator learning methods have been proposed to learn the mapping from the input space to the solution space, enabling rapid inference of solutions to partial differential equations (PDEs) once trained. In this work, we address the parametric elliptic interface problem. Building upon the deep operator network (DeepONet), we propose an extended interface deep operator network (XI-DeepONet). XI-DeepONet exhibits three unique features: (1) The interface geometry is incorporated into the neural network as an additional input, enabling the network to infer solutions for new interface geometries once trained; (2) The level set function associated with the interface geometry is treated as the input, on which the solution mapping is continuous and can be effectively approximated by the deep operator network; (3) The network can be trained without any input-output data pairs, thus completely avoiding the need for meshes of any kind, directly or indirectly. We conduct a comprehensive series of numerical experiments to demonstrate the accuracy and robustness of the proposed method.
		\end{abstract}
		
		\begin{keyword}
			Operator learning \sep Elliptic interface problem \sep Level set function
			
		\end{keyword}
	\end{frontmatter}


\section{Introduction}

Elliptic interface problems arise in numerous scientific domains, such as materials science \cite{greengard1994numerical} and fluid dynamics\cite{stroud2002numerical}, where the background comprises distinct materials on subdomains separated by smooth curves or surfaces, termed the interface. The solutions to interface problems may exhibit non-smoothness or even discontinuity due to the non-continuous coefficients and transmission conditions. The low global regularity of the solutions, coupled with the intricate geometry of the interfaces, poses significant challenges for numerical simulations.

Addressing these challenges, a myriad of numerical methods have been developed over the past several decades. These classical numerical methods for solving elliptic interface problems can be broadly categorized into two groups: fitted and unfitted methods. Fitted methods discretize interface problems using body-fitted meshes, ensuring that the interface does not intersect with mesh elements. Subsequently, the mesh can be used to apply the standard finite element method (FEM) \cite{zienkiewicz2005finite, mu2013weak}. Chen and Zou \cite{chen1998finite} conducted a study on the linear FEM on nearly fitted quasi-uniform meshes and proved optimal order error estimates up to some logarithm factors. Li et al. \cite{li2010optimal} extended this result to higher-order FEM. Guyomarc'h, Lee, and Jeon \cite{guyomarc2009discontinuous} analyzed the local discontinuous Galerkin Method. Burman and Hansbo \cite{burman2010interior} developed an interior penalty finite element method that employs a jump-stabilized Lagrange multiplier. Dryja considered two fitted linear discontinuous Galerkin methods for two-dimensional elliptic interface problems on piecewise fitted meshes obtained by independently partitioning each subdomain \cite{dryja2003discontinuous}. In \cite{he2020interface}, an interface penalty FEM (IPFEM) was proposed for elliptic interface problems, which allows the use of different meshes in different subdomains separated by the interface.

FEMs represent a viable approach for resolving interface problems on fitted meshes. Nevertheless, generating high-quality interface-fitted meshes for complex interfaces, particularly for time-dependent problems, is extremely costly. As a result, attempts have been made to solve elliptic interface problems using unfitted meshes. Unfitted methods allow the interface to cut through mesh elements, while special techniques are required to incorporate jump conditions across the interface with these methods. One approach is to use immersed finite element methods based on Cartesian meshes (see, e.g., \cite{leveque1994immersed, li1998immersed, li2003new, gong2008immersed, gong2010immersed, lin2019nonconforming}). In this method, the standard finite element basis functions are locally modified for elements intersected by the interface to accurately or approximately satisfy the jump conditions across the interface. The other approach is to use the extended finite element methods (XFEMs) based on unfitted-interface meshes, which are primarily designed to solve problems with discontinuities, kinks, and singularities within elements (see \cite{ hansbo2002unfitted, xiao2020high, zhang2022condensed} and references therein). For XFEMs, extra basis functions are added for elements intersected by the interface to capture discontinuous conditions, and the jump conditions are enforced using the variation of Nitsche's approach.

Recently, there has been a burgeoning interest in the scientific computing community in utilizing machine learning-based methods,
such as physics-informed neural networks (PINNs) \cite{raissi2019physics}, the deep Ritz method \cite{yu2018deep}, the random feature method \cite{chen2022rfm}, to solve interface problems. These methods apply neural network approximations to solve PDEs, offering the mesh-free advantage that traditional numerical methods may lack. They are particularly well-suited for handling complex interfaces, rendering them an ideal choice for solving elliptic interface problems. Lin et al. \cite{lai2022shallow} proposed a completely shallow Ritz network framework for solving interface problems by augmenting the level set function. The introduced shallow neural network with one hidden layer can significantly reduce the training cost in contrast to deep neural networks. Wu and Lu \cite{wu2022inn} introduced an interface neural network that decomposes the computational domain into several subdomains, with each network responsible for the solution in each subdomain. This method is effective in addressing the low global regularity of the solution. A similar piecewise deep neural network method was also introduced in \cite{he2022mesh}. A discontinuity capturing shallow neural network and a cusp-capturing PINN were developed in \cite{hu2022discontinuity, tseng2023cusp}. Both of them augment a coordinate variable to construct the loss function of interface conditions. The crucial idea is that a $d$-dimensional piecewise function can be extended to a continuous function defined in a $(d+1)$-dimensional space. The random feature method was applied to interface problems where two sets of random feature functions are employed to approximate the solution along each side of the interface \cite{chi2024random}.

The methods mentioned above are designed to solve specific interface problems with given conditions. If some conditions or the interface geometry change, one has to solve the interface problem again, leading to high computational costs and significant time consumption. Therefore, developing operator neural networks for elliptic interface problems is crucial, as they can provide rapid simulations for different input functions with a neural network trained once. Two notable examples are the deep operator network (DeepONet) \cite{lu2021learning} and Fourier neural operator (FNO) \cite{li2020fourier}. Both have been developed to directly learn the solution operator of arbitrary PDEs, mapping between two infinite-dimensional Banach spaces. Such methods hold significant potential for the development of rapid forward and inverse solvers for PDE problems and have demonstrated exceptional performance in numerous applications \cite{lu2022comprehensive, kontolati2024learning, tripura2023wavelet}. By leveraging automatic differentiation, physics-informed DeepONet (PI-DeepONet) \cite{wang2021learning} and physics-informed neural operator (PINO) \cite{li2024physics} have been proposed to learn the solution operator of PDEs, even in the absence of paired input-output training data. However, DeepONet and FNO are designed solely to learn operators defined on a single Banach space. To overcome this limitation, a novel operator learning network architecture called multiple-input operator network (MIONet) has been proposed for multiple input operators defined on the product of Banach spaces in \cite{jin2022mionet}. \cite{wu2024solving} introduced an interface operator network (IONet) to approximate the solution operator for parametric elliptic interface problems. IONet divides the entire domain into several subdomains and employs multiple branch networks and trunk networks. This method effectively captures discontinuities in both input functions and output solutions across the interface.

However, the piecewise neural network framework requires classifying the training points based on their respective subdomains and then inputting them into their corresponding branch and trunk networks. This implies that IONet can only solve a specific class of interface problems with a fixed interface. It is not capable of solving elliptic interface problems with varying interfaces, such as those involving positional or geometric changes. To overcome this limitation, we utilize a continuous level set function to effectively represent the position of the interface and develop a new DeepONet architecture capable of solving elliptic interface problems with parameterizable interfaces. In both classical numerical methods \cite{gong2008immersed} and PINNs \cite{tseng2023cusp, lai2022shallow}, the level set function plays a crucial role in solving elliptic interface problems. Inspired by the approach of augmenting the level set function as an additional feature input in \cite{tseng2023cusp}, we reformulate the equation to establish a connection with the level set function. We then use the level set function, which represents the geometric shape of the interface, as the input function for a branch network. This allows the neural network to express the potential relationship between the interface position information and the solutions, enabling it to handle evolving interfaces. Additionally, we propose a physics-informed loss to significantly reduce the need for extensive training data, allowing the network to perform effectively without requiring paired input-output observations within the computational domain. Once trained, this method enables quick simulations with various input functions and interfaces. We would like to emphasize that our methods can be implemented within any operator learning framework, such as FNO, to achieve similar results. However, we have chosen to use the DeepONet framework in this particular study due to its ability to handle problems involving complex interfaces.

The rest of this paper is organized as follows. In Section \ref{Section 2}, we introduce the elliptic interface problem and provide an overview of the fundamental concepts underlying DeepONet. In Section \ref{Section 3}, we present a novel neural network designed to address elliptic problems with varying interfaces. We will demonstrate the integration of geometric information of the interface into the neural network by augmenting the level set function in the trunk network and setting the level set function as the input function in the branch network, respectively. Numerical results are presented in Section \ref{Section 4} to illustrate the efficiency of the proposed method. Finally, we give concluding remarks in Section \ref{Section 5}.

\section{Elliptic interface problem and operator neural networks}
\label{Section 2}
In this section, we first present the elliptic interface problem, followed by a brief introduction to the MIONet and PI-DeepONet. Both MIONet and PI-DeepONet are extensions of DeepONet, and they will play significant roles in our work.

\subsection{Interface problem}
In this paper, we will consider the numerical solution of the elliptic interface problems with discontinuous coefficients by the deep learning methods. Let $\Om=\Om^+\cup\Ga\cup\Om^-$ be a bounded domain in $\R^d, d\ge 1$, where $\Om^+$ and $\Om^-$ are {two subdomains} of $\Om$ and $\Ga=\pa\Om^+\cap\pa\Om^-$. 
Consider the following elliptic interface problem:
\begin{align}
- \nabla  \cdot ( {a(x) \nabla u} ) + b(x) u &= f \qquad \mathrm{in}\;\Omega \backslash \Gamma, \label{2.1}\\
{\left[ u \right]_\Gamma } &= g_D \quad\,\, \mathrm{on}\;\Gamma, \label{2.2}\\
{\left[ {a( x ){\partial _{\bn}}u} \right]_\Gamma } &= {g_N} \quad\,\, \mathrm{on}\;\Gamma, \label{2.3}\\
u &= h \qquad \mathrm{on}\; \partial \Omega. \label{2.4}
\end{align}
Here, we have used $\partial_{\bn} u$ to denote the normal derivative $\nabla u \cdot {\bn}$, where ${\bn}$ is the unit normal vector to $\Gamma$ from $\Omega^-$ to $\Omega^+$. The functions $a,f:\Omega \to \mathbb{R}$, $ g_N $, $g_D:\Gamma \to \mathbb{R},$ and $h:\partial \Omega \to \mathbb{R}$ are assumed to be bounded {and $a(x) \ge a_0 > 0$ for all $x \in \Omega$.} We denote the restrictions of $\xi$ on $\Omega^+$ and $\Omega^-$ for any function $\xi:\Omega \to \mathbb{R}$ by
$$
{\xi ^ + } = {\left. \xi  \right|_{{\Omega ^ + }}} \quad \mathrm{and} \quad {\xi ^ - } = {\left. \xi  \right|_{{\Omega ^ - }}}
$$
respectively, and we denote by
$$
{\left[ \xi  \right]_\Gamma }( {x} ) = \mathop {\lim }\limits_{y \to x,y \in {\Omega ^ + }} {\xi ^ + }( y ) - \mathop {\lim }\limits_{y \to x,y \in {\Omega ^ - }} {\xi ^ - }( y )
$$
the jump of $\xi$ across the interface $\Gamma$, when the unique limiting values of $\xi$ from both sides of $\Gamma$ exists.
\subsection{Learning operators via neural networks}

MIONet proposed in \cite{jin2022mionet} is a method to learn the (typically nonlinear) operator mapping from a product of $n$ infinite dimensional Banach spaces (input spaces) to another infinite dimensional Banach space (output space). Let $\mathcal{V}_1, \mathcal{V}_2,...,\mathcal{V}_n$ be $n$ different input Banach space and $\mathcal{U}$ be a output space. We use $\mathcal{G}$ to denote the operator that maps between input spaces and output space. Then MIONet aims to learn a continue operator
$$
\mathcal{G}:{\mathcal{V}_1} \times {\mathcal{V}_2} \times  \cdots {\mathcal{V}_n} \to \mathcal{U}, \qquad( v_1, v_2,..., v_n) \mapsto  u,
$$
where $v_i \in \mathcal{V}_i$ and $u=\mathcal{G}( v_1, v_2,..., v_n) \in \mathcal{U}$. Based on the main neural network approximation theorems (see \cite[Theorem 2.5 and Corollary 2.6]{jin2022mionet}), we can construct MIONet $\mathcal{G}_\theta$ to approximate the operator $\mathcal{G}$.
\begin{align}
\mathcal{G}_\theta( {{v_1},...,{v_n}} )( x ) &= \mathcal{S}( {\underbrace {{N_{{b_1}}}( {{v_1}( y )} )}_{Branc{h_1}} \odot  \cdots  \odot \underbrace {{N_{{b_n}}}( {{v_n}( y )} )}_{Branch_n} \odot \underbrace {{N_t}( x )}_{Truck}} ) + \underbrace {{b_0}}_{bias}
\nonumber \\
&= \sum\limits_{i = 1}^m {{t^{(i)}}} \prod\limits_{j = 1}^n {b_j^{( i )} + {b_0}}, \nonumber
\end{align}
where $ \odot $ is the Hadamard product, and $\mathcal{S}$ denotes the summation of all the components of a vector.  $\theta$ represents the parameters of the neural network, i.e., all trainable parameters of sub-branch networks $\mathcal{N}_{b_i}, i=1,2,...,n,$ and trunk network $\mathcal{N}_t$ and the trainable bias $b_0 \in \mathbb{R}$.  $\left[ b^{(1)}_i, b^{(2)}_i,..., b^{(m)}_i \right]^T$ and $\left[ t^{(1)}, t^{(2)},..., t^{(m)}\right]^T $ represent the output of sub-branch network $\mathcal{N}_{b_i}$ and trunk network $\mathcal{N}_t$, respectively. In order to input the function $( v_1, v_2,..., v_n) \in ( {\mathcal{V}_1} \times {\mathcal{V}_2} \times  \cdots {\mathcal{V}_n}) $ into neural network, we need to discretize the input function by a collection of $k$ fixed point locations $\left\lbrace {y_1}, {y_2},...,{y_k}\right\rbrace $ called ``sensors''. Each input function dataset $v_i( {y}) := \left[ v_i( {y_1}), v_i( {y_2}),..., v_i( {y_k})\right] ^T$ is finite dimensional. 

As the usual DeepONet, the MIONet method requires a large corpus of paired input-output observations to train neural networks by minimizing the following empirical loss function.
$$
Loss(\theta) = \dfrac{1}{NP}{\sum\limits_{i = 1}^N {\sum\limits_{j = 1}^P {\left| {{\mathcal{G}_\theta }\left(  {{v_1}^{( i )},...,{v_n}^{( i )}} \right) \left( {{x}_j^{( i )}} \right)  - \mathcal{G}\left(  {{v_1}^{( i )},...,{v_n}^{( i )}} \right) \left(  {{x}_j^{( i )}} \right) } \right|^2} }},
$$
where $\left\lbrace {{v_1}^{( i )},...,{v_n}^{( i )}} \right\rbrace_{i=1}^N$ denotes $N$ input function groups sampled from the parameter space ${\mathcal{V}_1} \times {\mathcal{V}_2} \times  \cdots {\mathcal{V}_n}$, and $\left\lbrace x_j^{(i)} \right\rbrace _{j=1}^P$ are $P$ training data points for corresponding input function group $\left( {v_1}^{( i )},...,{v_n}^{( i )} \right) $. 

However, generating sufficient large training datasets may be relatively expensive and challenging. In \cite{wang2021learning}, the authors developed a method called PI-DeepONet, which is able to train such model without any observed data at all. It only depends on the given knowledge of the PDEs form and its corresponding initial and boundary conditions (IBCs). To explain this, we consider general parametric PDEs taking the form
$$
\mathcal{L}\left( {v_1, v_2,..., v_n, u} \right) = 0.
$$
Then the trainable parameter $\theta$ of neural networks can be optimized by minimizing the residuals of the equation using the automatic differentiation \cite{baydin2018automatic} like PINNs. Consequently, we may then construct a PI-DeepONet  by formulating the following loss function
$$
{L_{physics}}( \theta  ) = \dfrac{1}{{NP}}{\sum\limits_{i = 1}^N {\sum\limits_{j = 1}^P {\left| {\mathcal{L}\left(  {  {v_1^{( i )},v_2^{( i )},...,v_n^{( i )}}, \mathcal{G}\left( {v_1^{( i )},v_2^{( i )},...,v_n^{( i )}} \right) } \right) \left(  {{x}_j^{( i )}} \right) } \right|^2} }}.
$$

\section{Extended interface DeepONet method}
\label{Section 3}
In this section, we first show how to extend the solution function by augmenting the level set function, and then propose the neural network architecture of XI-DeepONet and its loss functions. For our method, it needs to assume that the interface condition~\eqref{2.2} is homogenous, i.e. $g_D=0$ ( See Remark~\ref{remark1} below for the reason). 

\subsection{Extension function by augmenting the level set function}
Now, we review the level set function of the interface $\Gamma$ which is assumed to belong to $C^{1,1}$. Let $\phi: \Omega \to \mathbb{R}$ be a smooth function that satisfies
\begin{align*}
\phi\left( x\right) =\left\{
\begin{matrix}
\begin{aligned}
&< 0,\\
&=0,\\
& >0,
\end{aligned}
&
\begin{aligned}
&\mathrm{if}\;x \in \Omega ^ - ,\\
&\mathrm{if}\;x \in \Gamma ,\\
&\mathrm{if}\;x \in \Omega ^ + .
\end{aligned}
\end{matrix}
\right.
\end{align*}

To homogeneous jump condition problem, let $U(x,z)\in \mathbb{R}$ be a continuous function defined on $\Omega \times \mathbb{R}$. Assume that $U({{x},z=\Phi({x})})= u(x)$ with \eq{\label{phifun}
	\Phi ( {x}) := \left| {\phi ( {x} )} \right|.}
That is to say, the domain of $u$ can be regarded as the projection of the domain of $U$. Thus, the gradient of $u$ can be calculated as follows:
\begin{equation}
\label{eq3.1} \nabla u = {\nabla _{{x}}}U + {\partial _z}U\nabla \Phi,
\end{equation}
where ${\nabla _{{x}}}U$ represents a vector with partial derivatives of $U$ with respect to the components in ${x}$ and $\partial_z U$ is the partial derivative of $U$ with respect to $z$. It is easy to see
\begin{align*}
\nabla \Phi =\left\{
\begin{matrix}
\begin{aligned}
&\nabla \phi ,\\
&-\nabla \phi ,\\
\end{aligned}
&
\begin{aligned}
&\mathrm{in}\;{\Omega ^ + },\\
&\mathrm{in}\;{\Omega ^ - },
\end{aligned}
\end{matrix}
\right.
\end{align*}
Furtherly, using (\ref{eq3.1}) we can get the following PDEs with homogeneous jump condition,
\begin{equation}
\mathcal{L}\left[ u \right]  = 0  \quad \mathrm{in}\;\Omega \backslash \Gamma, \qquad \mathcal{I}\left[ u \right]=0\quad \mathrm{on}\; \Gamma, \qquad \mathcal{B}\left[ u \right] = 0 \quad \mathrm{on}\; \partial \Omega,
\nonumber
\end{equation}
where
\begin{align}
\mathcal{L}\left[ u \right] := &- \nabla  \cdot ( {a \nabla u} ) + bu - f \nonumber \\
= &-a( {{\Delta _{{x}}}U + 2\nabla \Phi  \cdot {\nabla _{{x}}}( {{\partial _z}\Phi } ) + {{\left| {\nabla \Phi } \right|}^2}{\partial _{zz}}U + {\partial _z}U\Delta \Phi } ) \label{3.3}\\
&- \nabla a \cdot ( {{\nabla _{{x}}}U + {\partial _z}U\nabla \Phi } ) + bU - f, \nonumber
\end{align}
and
\begin{align}
\mathcal{I}[ u ] &:= {\left[ {a( x ){\partial _{\bn}}u} \right]_\Gamma } -  {g_N}
= {\left[ {( {a{\nabla _{{x}}}U + {\partial _z}\nabla \Phi } ) \cdot n} \right]_\Gamma } -  {g_N},  \label{3.4}\\
\mathcal{B}[ u ] &:= u-h = U-h \label{3.5}. 
\end{align} 
Note here $\Delta_{{x}}$ is the Laplace operator concerning only the variable ${x}$. 

\begin{remark}
	\label{remark1}
	Since the level set function is continuous throughout the entire computational domain $\Omega$, the functions represented by neural networks must inherently be continuous at the interface. This rationale leads to our assumption that the jump condition in \eqref{2.2} is homogeneous. For elliptic interface problems with non-homogeneous jump condition $\left[ u \right]_\Gamma =g_D$, we can extend the jump function $g_D:\Gamma \to \mathbb{R}$ to a piecewise smooth function $v:\Omega \to \mathbb{R}$ such that
	$$
	{\left[ v \right]_\Gamma } = {g_D} \qquad \mathrm{and} \qquad v \equiv 0 \quad in\; \Omega^-.
	$$
	Thus, the function $w:=u-v$ satisfies
	an interface problem with homogeneous jump condition, hence it can be solved by our proposed method.
\end{remark} 

\subsection{Neural network architecture of XI-DeepONet and its function}
In this subsection, we present an operator learning neural network to approximate the solution $U({x}, \Phi({x}))$. The core idea of our novel method involves utilizing the level set function providing the geometric feature of the interface as the input function. Consequently, this approach enables us to obtain numerical solutions for elliptic interface problems with different interfaces after training once. For an illustration, we consider the operator $\mathcal{G}$ mapping from source term $f$ and the absolute value function of level set function $\Phi$ to the solution $u$,
$$
\mathcal{G}:( f( {x} ) , \Phi ( {x} ))  \to u(x).
$$
We define the following XI-DeepONet to approximate the operator $\mathcal{G}$:
$${\mathcal{G}_\theta }(f,\Phi)(x) := {G_\theta }( {f,\Phi } )( {{x},\Phi ( {x} )} ) 
= \left\{ {\begin{array}{*{20}{lc}}
	{{G_\theta }( {f,\Phi } )( {{x},\phi ( {x} )} ),}&{\mathrm{if}\;{x} \in {\Omega ^ + }},\\
	{{G_\theta }( {f,\Phi } )( {{x}, - \phi ( {x} )} ),}&{\mathrm{if}\;{x} \in {\Omega ^ - }}.
	\end{array}} \right.$$
Similarly, following the MIONet method, the operator is constructed via Hadamard product, summation, and a bias. Its formula is defined as follows:
\begin{align}
\mathcal{G_{\theta}}( {f,\Phi } )( {x} ){\rm{ }} &= \mathcal{S}( {\underbrace {{N_{{b_1}}}( {{f}( {y} )} )}_{Branch_1} \odot \underbrace {{N_{{b_2}}}( {{\Phi}( {y} )} )}_{Branch_2} \odot \underbrace {{N_t}( {x} )}_{Truck}} ) + \underbrace {{b_0}}_{bias}
\nonumber \\
&= \sum\limits_{i = 1}^m {{t^{(i)}}} \prod\limits_{j = 1}^{2} {b_j^{( i )} + {b_0}}. \nonumber
\end{align}

\begin{figure}[htbp]
	\centering
	\includegraphics[scale=0.35]{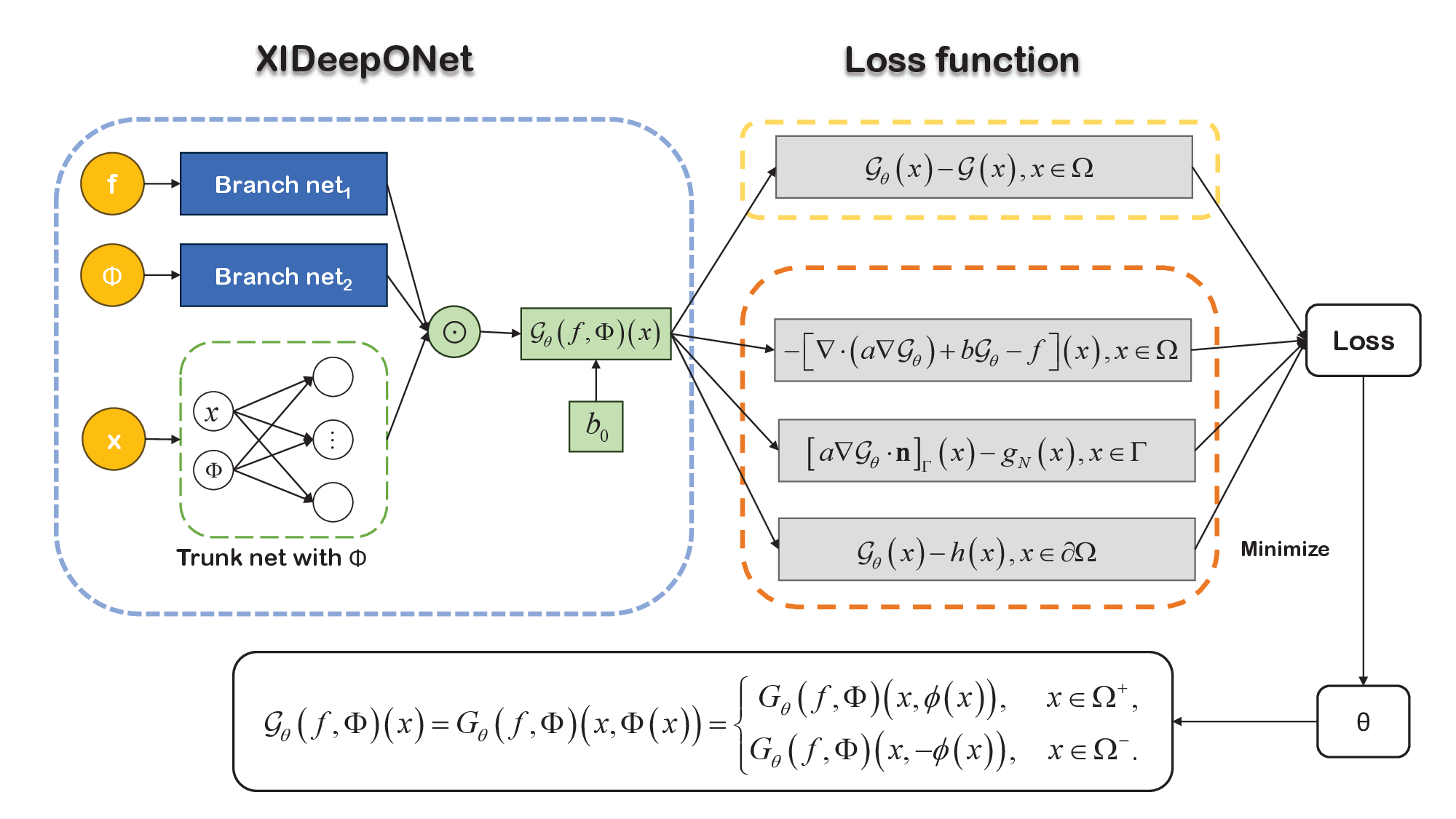}
	\caption{The network architecture of Extended Interface DeepONet.}
	\label{Neunal_network}
\end{figure}

In this paper, we only use the simplest feedforward neural network (FNN) to demonstrate the capability of XI-DeepONet for solving elliptic interface problems with varying interfaces. Similar to MIONet, a neural network approach can be adopted to optimize the parameters $\theta$ by minimizing the following loss function:
\begin{equation}\label{lossdd}
L_{operator}(\theta) = \dfrac{1}{{NP}}{\sum\limits_{i = 1}^N {\sum\limits_{j = 1}^P {\left| {{\mathcal{G}_\theta }( {{{f}^{(i)}},{{\Phi} ^{(i)}}} )( {{x}_j^{(i)}} )-\mathcal{G}( {{{f}^{(i)}},{{\Phi}^{(i)}}} )( {{x}_j^{(i)}} )} \right|^2}}}.
\end{equation}
A data-driven (DD) approach can be adopted to train the network by $L_{operator}(\theta)$. However, the cost of obtaining a large amount of paired experimental data and high-precision numerical simulation is generally expensive. Therefore, we will use the physics information neural network which optimizes the parameters $\theta$ by minimizing the following composite loss function
\begin{equation}\label{lossph}
{L_{physics}}( \theta  ) =  {\beta _{{\Omega}}}{L_{{\Om}}}( \theta  ) + {\beta _{\partial \Omega }}{L_{\partial \Omega }}( \theta  ) + {\beta _\Gamma }{L_\Gamma }( \theta  ),
\end{equation}
where the residual error $L_{\Om}$, boundary condition error $L_{\partial\Om}$ and interface condition error $L_\Gamma$ are defined as follows:
\begin{equation*}
{L_{{\Om }}}( \theta  )  = {\sum\limits_{i = 1}^N {\sum\limits_{j = 1}^{{P_\Om}} {\left| \mathcal{L}\left[ {{\mathcal{G}_\theta }( {{f^{( i )}},{\Phi ^{( i )}}} )( {{x}_{j,\Om}^{( i )}} )}  \right]  \right| ^2} }},  
\end{equation*}
\begin{equation*}
{L_{\pa\Om}}( \theta  ) = {\sum\limits_{i = 1}^N {\sum\limits_{j = 1}^{P_{\pa\Om}} {\left| \mathcal{B}\left[ {{\mathcal{G}_\theta }( {{f^{( i )}},{\Phi ^{( i )}}} )( {{x}_{j,\pa\Om}^{( i )}} )}  \right] \right| ^2} }},
\end{equation*}
\begin{equation*}
{L_\Gamma }( \theta  ) = {\sum\limits_{i = 1}^N {\sum\limits_{j = 1}^{{P_\Gamma}} {\left|  \mathcal{I}\left[ {{\mathcal{G}_\theta }( {{f^{( i )}},{\Phi ^{( i )}}} )( {{x}_{j,\Gamma}^{( i )}} )}  \right] \right| ^2}}}.
\end{equation*}
Here, for the i-th input function, $\left\lbrace {x}_{j,\Om}^{(i)} \right\rbrace_{j=1}^{P_\Om}$ are the domain training points sampled from the computational domain $\Omega$, $\left\lbrace {x}_{j,\pa\Om}^{(i)} \right\rbrace_{j=1}^{P_{\pa\Om}}$ are randomly sampled from the boundary $\partial \Omega$, and $\left\lbrace {x}_{j,\Gamma}^{(i)} \right\rbrace_{j=1}^{P_\Gamma}$ represent the training data points randomly sampled from the interface $\Gamma$. The constants $\be_\Om, \be_{\pa\Om}$ and $\be_\Ga$ appeared in the loss function (\ref{lossph}) are chosen to balance the contribution of the terms related to the boundary condition and  interface jump condition, respectively. Fig.~\ref{Neunal_network} presents a schematic visualization of the proposed network architecture of XI-DeepONet. 
\begin{remark}
	\label{remark2}
	It is important to note that the choice of $\Phi$ is not unique, which depends on the level set function. We can avoid the singularity of the level set function $\phi$ by defining properly the function $\Phi$. For instance (see Example~\ref{exa3}), if the derivative of $\phi(x)$ is singular at certain points in $\Omega^-$, we can use the ReLU function instead of the absolute value function, i.e., $\Phi({x}):=ReLU(\phi({x}))$. Therefore, $\Phi(x)$ and its derivative are equal to $0$ if $x \in \Omega^-$. 
\end{remark}

\section{Numerical results}
\label{Section 4}
In this section, we aim to demonstrate the capability of XI-DeepONet for solving the elliptic interface problems, i.e., equations (\ref{2.1})--(\ref{2.4}), with different interfaces. For all numerical experiments,  both branch networks and trunk networks are FNNs. All operator network models are trained via stochastic gradient descent using Adam optimizer \cite{kingma2014adam}, and the learning rate is set to exponential decay. Based on the experience of solving elliptical interface problems using PINNs, we set hyperparameters $\beta_{{\Omega }}=\beta_\Gamma= 1$ and $\beta_{\partial \Omega }=100$, respectively.  Furthermore, XI-DeepONets trained by minimizing the loss function (\ref{lossdd}) and (\ref{lossph}) are denoted by ``DD-XI-DeepONet'' and ``PI-XI-DeepONet'' respectively. The training data points set $M$ is randomly sampled in the computational domain, and test data points set $M_{test}$ comes from equispaced grid points. We calculate the relative $L^2$ errors defined by  ${{{{\left\| {{u_S} - u} \right\|}_{2} }} \mathord{\left/{\vphantom {{{{\left\| {{u_S} - u} \right\|}_{2} }} {{{\left\| u \right\|}_{2} }}}} \right.\kern-\nulldelimiterspace} {{{\left\| u \right\|}_{2} }}}$, where
$$
{\left\| u \right\|_2} = \sqrt {\dfrac{1}{{{M_{test}}}}\sum\limits_{i = 1}^{{M_{test}}} {{{( {u( {{{x}_i}} )} )}^2}}},
$$
and $u_S$ denotes the numerical solution obtained by the neural network. In general, we set the maximum iteration step $epoch=4\times10^{4}$. Branch and trunk networks are 5-layer FNNs with 100 units per layer with activation function $ReLU$ and $Tanh$ for DD-XI-DeepONet and PI-XI-DeepONet, respectively. All trials are run on an NVIDIA RTX4090 GPU.

\begin{figure}[htbp]
	\centering
	\begin{minipage}{0.45\linewidth}
		\centering
		\includegraphics[width=1\linewidth]{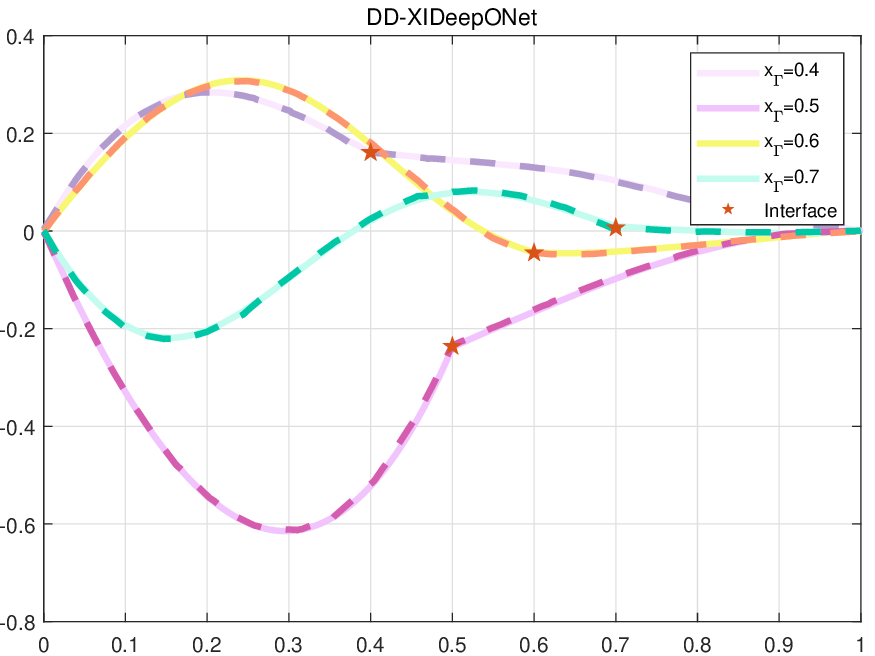}
	\end{minipage}
	\begin{minipage}{0.45\linewidth}
		\centering
		\includegraphics[width=1\linewidth]{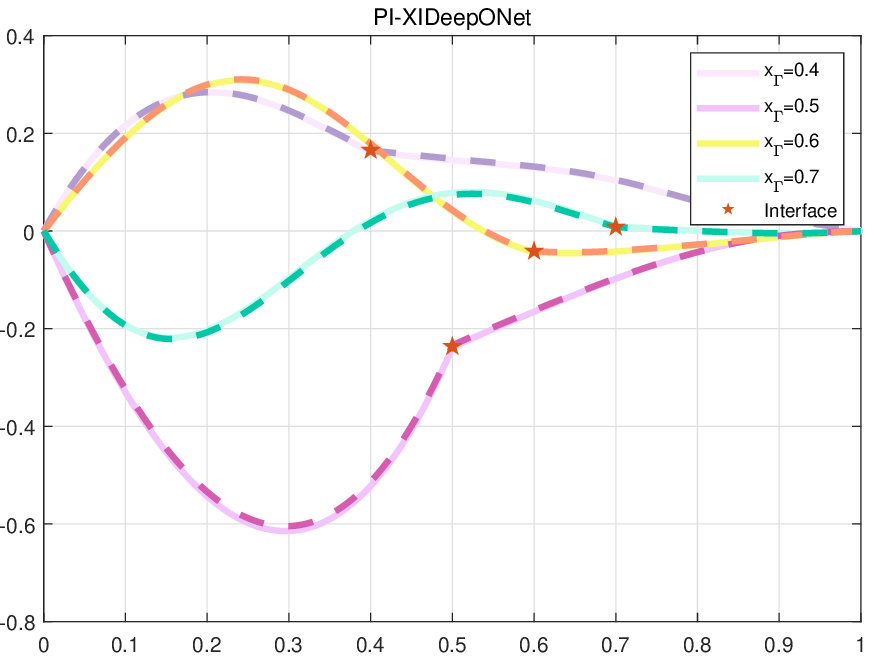}
	\end{minipage}
	\caption{Four representative curve of solution from the test sample data (distinguished by different colors). The reference solutions obtained by MIB and the predicted solutions computed by NNs are depicted by solid and dashed lines respectively. Left: The approximate solution acquired by DD-XI-DeepONets. Right: The approximate solutions acquired by PI-XI-DeepONets.}
	\label{Fig2}
\end{figure}

\begin{example}\label{exa1}
	Consider the following one-dimensional Possion equation on the interval $\Omega=[0,1]$ with a random interface point $x_\Gamma=p\in [0.4, 0.7]$:
	\begin{align*}
	- (& {a(x) u_x( x )} )_x = f( x ), \quad x\in [0,1], \\
	&\jm{u}\big|_{x_\Ga}=0, \quad  \jm{a(x) u_x}\big|_{x_\Ga}=0, \\
	&	u( 0) =u( 1) =0,
	\end{align*}
	where the coefficient $a(x)$ is a piecewise constant defined by
	$$
	a( x ) = \left\{ {\begin{array}{*{20}{lc}}
		{0.1,}&{\quad \mathrm{in}\; \Omega^-:=[0,x_\Ga),}\\
		{0.5,}&{\quad \mathrm{in}\; \Omega^+:=(x_\Ga,1].}
		\end{array}} \right.
	$$
	and the level set function is denoted by
	$$
	\phi\left( x\right) = x - p.
	$$
\end{example}
This example is designed to demonstrate the capability of the present network XI-DeepONet. The solution of this example is non-smooth across the interface point, and our objective is to learn the solution operator $\mathcal{G}$ mapping from source function $f$ and level set function $\phi$ to the latent solution $u(x)$. The value of $p$ is a random real number sampled from the interval $[0.4, 0.7]$. Further, we model the input function $f$ using a zero-mean Gaussian random field (GRF) \cite{lu2021learning} with the following covariance kernel
$$
{k_l}( {{x_1},{x_2}} ) = \exp ( { - {{{{\left\| {{x_1} - {x_2}} \right\|}^2}} \mathord{\left/
			{\vphantom {{{{\left\| {{x_1} - {x_2}} \right\|}^2}} {2{l^2}}}} \right.
			\kern-\nulldelimiterspace} {2{l^2}}}} ).
$$
We use two GRFs with length scales $l_1=0.2$ and $l_2=0.1$ to independently simulate $f^+$ and $f^-$, respectively (cf. \cite{wu2024solving}). Since the analytical solution is unknown, the ``exact'' solutions for training and testing come from the numerical solution via the matched interface and boundary (MIB) method \cite{yu2007matched} on the $1 \times 1000$ grid. We randomly sample $10000$ and $1000$ input functions $f$ and $\Phi$ for training and testing, respectively. The $100$ locations of sensors are uniformly distributed in the interval $[0, 1]$. Note that paired input-output measurements of each training sample used to train DD models are generated by solving the given equation via the MIB method also. The information on neural networks and the mean relative $L^2$ errors between the approximate and reference solutions are presented in Table \ref{Table1}. Furthermore, we depict the approximate solutions of neural networks for four locations of the interface ($x_{\Gamma}=0.4$, $x_{\Gamma}=0.5$, $x_{\Gamma}=0.6$, $x_{\Gamma}=0.7$) in Fig.~\ref{Fig2}. 
\begin{table}[htbp]
	\centering
	\begin{tabular}{c|cccc}
		\toprule
		Model     & Activation & 
		Sensors   &
		Parameters & $\dfrac{\left\| {{u_S} - u} \right\|_{2}}{\left\| u \right\|_{2}}$\\
		\midrule
		DD-XI-DeepONet & ReLU &    100  & 141701   & $1.72 \times 10^{-2}$\\
		PI-XI-DeepONet & Tanh  &    100   & 141701   & $1.37 \times 10^{-2}$\\
		\bottomrule
	\end{tabular}%
	\caption{Example~\ref{exa1}: The mean of relative $L^2$ error of 1000 test function for DD and PI XI-DeepONet. }
	\label{Table1}
\end{table}

\begin{figure}[htbp]
	\centering
	\begin{minipage}{0.45\linewidth}
		\centering
		\includegraphics[width=1\linewidth]{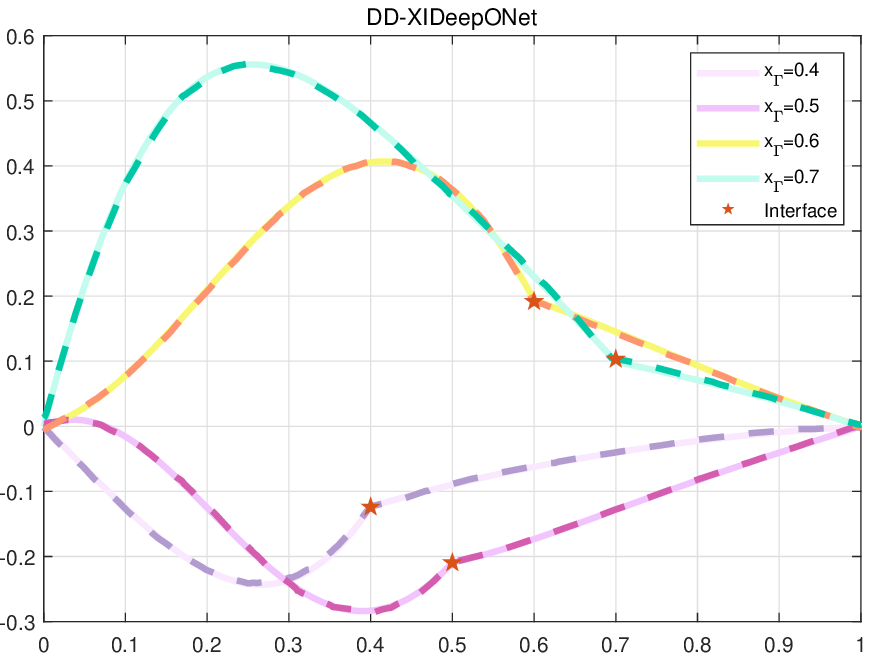}
	\end{minipage}
	\begin{minipage}{0.45\linewidth}
		\centering
		\includegraphics[width=1\linewidth]{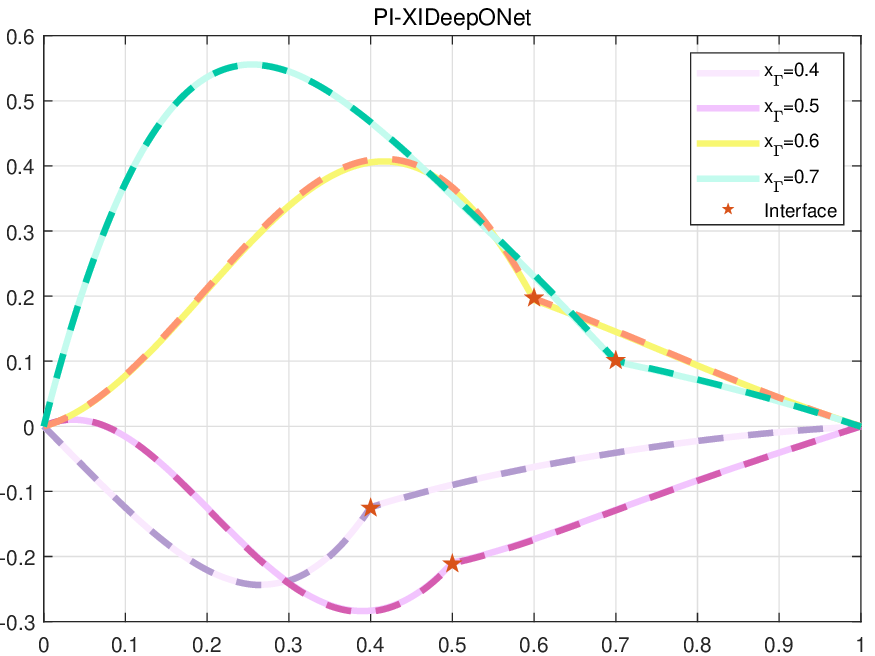}
	\end{minipage}
	\begin{minipage}{0.45\linewidth}
		\centering
		\includegraphics[width=1\linewidth]{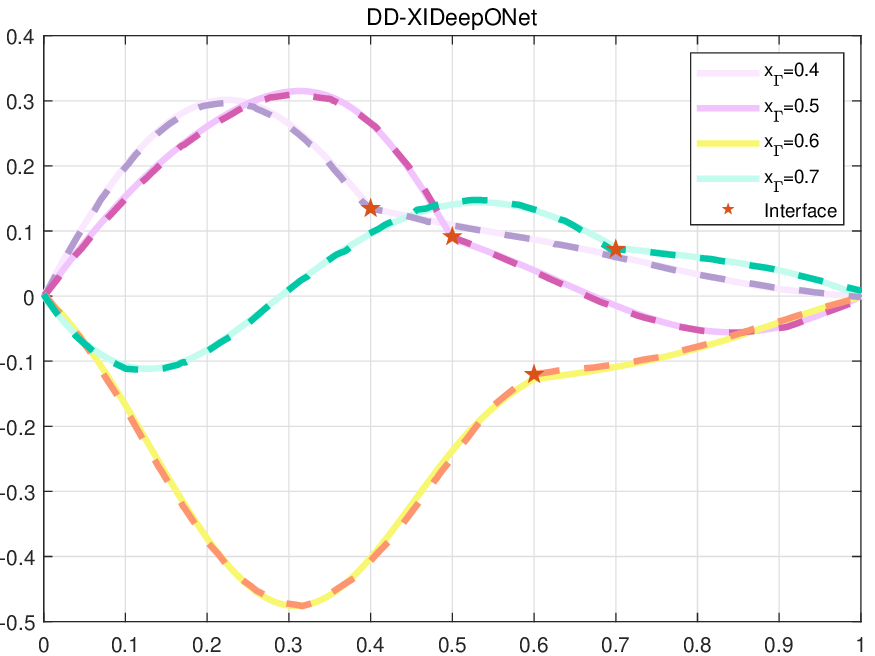}
	\end{minipage}
	\begin{minipage}{0.45\linewidth}
		\centering
		\includegraphics[width=1\linewidth]{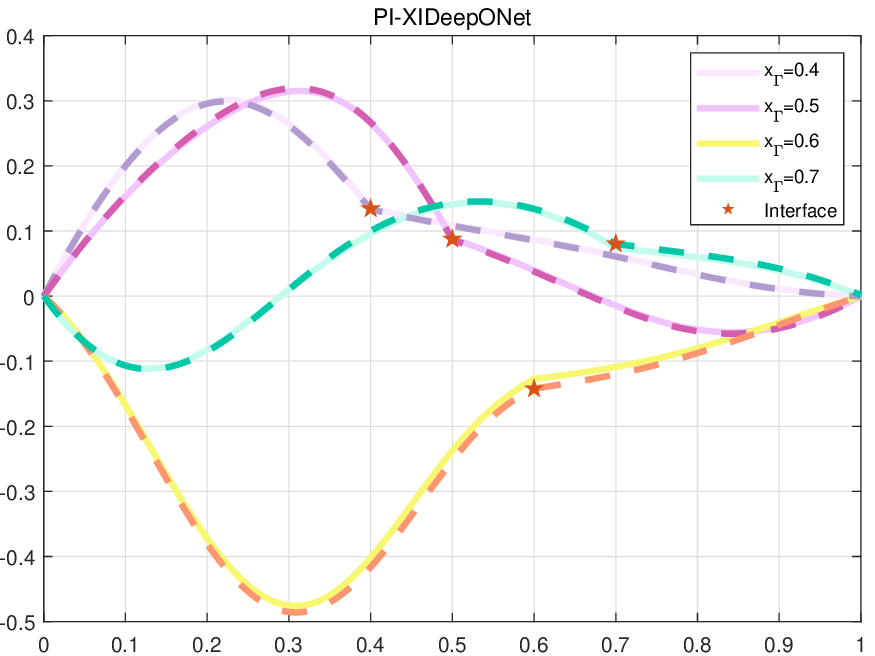}
	\end{minipage}
	\caption{Top: $l_1=l_2=0.3$. Bottom: $l_1=l_2=0.15$. The solid and dashed lines are the reference solution and the predicted solutions obtained by MIB and NNs, and the training data is generated by $l_1=0.2, l_2=0.1$.}
	\label{Fig3}
\end{figure}

It is observed that XI-DeepONet is highly effective in obtaining approximate solutions. When we use the test source function generated by $l_1=l_2=0.3$ (first row) and $l_1=l_2=0.15$ (second row), the model is still valid as shown in Fig.~\ref{Fig3}. Note that the training input source functions are generated only by $l_1=0.2$ and $l_1=0.1$. These results highlight the robustness of XI-DeepONet.

\label{Comparison}
In addition, we compare the performance of IONet \cite{wu2024solving}, DeepONet and XI-DeepONet under almost the same settings through this example. We still employ two GRFs with length scale $l_1=0.2$ and $l_2=0.1$ to simulate source terms and randomly sample 10000 input functions for training and 1000 input functions for testing, respectively. All neural network frameworks are trained by minimizing the residuals of interface equations and their corresponding boundary conditions. We fix the interface position $x_\Gamma = 0.5$ and only use the source term as input function without regard to variation of interface position. Table \ref{Table2} displays the different network architectures used to approximate the solution operator and records errors of the three models after 40K gradient descent iterations. It can be observed that XI-DeepONet, similar to IONet, demonstrates a strong capability to capture the unsmooth property of real solutions at the interface.
\begin{table}[htbp]
	\centering
	\begin{tabular}{c|ccccc}
		\toprule
		Model     & Activation & 
		Depth   &
		Width & No. of parameters &$\dfrac{\left\| {{u_S} - u} \right\|_{2}}{\left\| u \right\|_{2}}$\\
		\midrule
		XI-DeepONet & Tanh  &    5  & 100  & 141701 & $3.3 \times 10^{-3}$\\
		IONet & Tanh  &    5   & 90 & 140582 & $2.7 \times 10^{-3}$\\
		DeepONet & Tanh  &    5   & 140 & 149631 & $3.5 \times 10^{-1}$\\
		\bottomrule
	\end{tabular}%
	\caption{Example~\ref{exa1}: The mean of relative $L^2$ error of 1000 test function for three problems. }
	\label{Table2}
\end{table}

However, XI-DeepONet has the ability to generalize interface positions. This means that we can handle multiple problems with different interface. For IONet, the limitations of its framework preclude the implementation of this feature. To elucidate this limitation, we uniformly select 9 points within the interval $\left[0.3, 0.7\right] $ as interface points and train the model with 1000 input functions at each interface point. That is to say, for IONet, it requires to be trained nine times. As shown in the Table \ref{Table3}, XI-DeepONet saved nearly half of the time while achieving similar accuracy.
\begin{table}[htbp]
	\centering
	\begin{tabular}{c|ccccc}
		\toprule
		Model   & 
		Depth   &
		Width & No. of parameters & Time(s) &$\dfrac{\left\| {{u_S} - u} \right\|_{2}}{\left\| u \right\|_{2}}$\\
		\midrule
		XI-DeepONet  &    5  & 100  & 141701 & 4856 & $5.2 \times 10^{-3}$\\
		IONet  &    5   & 90 & 140582 & 1273 $\times$ 9 & $3.9 \times 10^{-3}$\\
		\bottomrule
	\end{tabular}%
	\caption{Example~\ref{exa1}: The mean of relative $L^2$ error and computational cost in seconds (s) for trained models. }
	\label{Table3}
\end{table}

\begin{example}\label{exa2}
	Consider the model problem \eqref{2.1}--\eqref{2.4} in the circular domain with a radius $R=1$. The embedded interface $\Gamma$ is a circle described by the zero level set of the function $\phi(x,y)=x^2 + y^2 - r_0^2$, where $r_0$, the radius of interface, belongs to $[0.4, 0.8]$. Further, the diffusion coefficients $a^+=1000, a^-=1$ are defined in a piecewise constant manner. The interface conditions on $\Gamma$ are given as
	$$
	g_D(x,y)=0, \quad g_N(x,y)=0,
	$$
	and $h(x,y)=0$ on $\partial \Omega$. 
\end{example}	

This example is designed to illustrate the applicability of the proposed method to solving high contrast interface problems. Here, the level set function and the source term $f$ are the input parameters of the target solution operator. In this example, we model the input function in the following way:
$$
f( {x,y} ) = 9p\sqrt {{x^2} + {y^2}} ,
$$
where $p \in [0.5, 2.0]$. When $p=1$, we obtain the exact solution $u$ as follows:
\begin{align*}
u\left( x, y\right) =\left\{
\begin{aligned}
&\dfrac{1}{a^+}\left(x^2+y^2\right)^{3/2}-\dfrac{1}{a^+ } &\mathrm{in}\;\Omega^+,\\
&\dfrac{1}{a^-}\left(x^2+y^2\right)^{3/2}+ \left( \dfrac{1}{a^+}-\dfrac{1}{a^-}\right)r_0^3-\dfrac{1}{a^+ } &\mathrm{in}\;\Omega^-,
\end{aligned}
\right.
\end{align*}
which will be used as the ground truth for later testing.

\begin{figure}[htbp]
	\centering
	\includegraphics[width=0.4\linewidth]{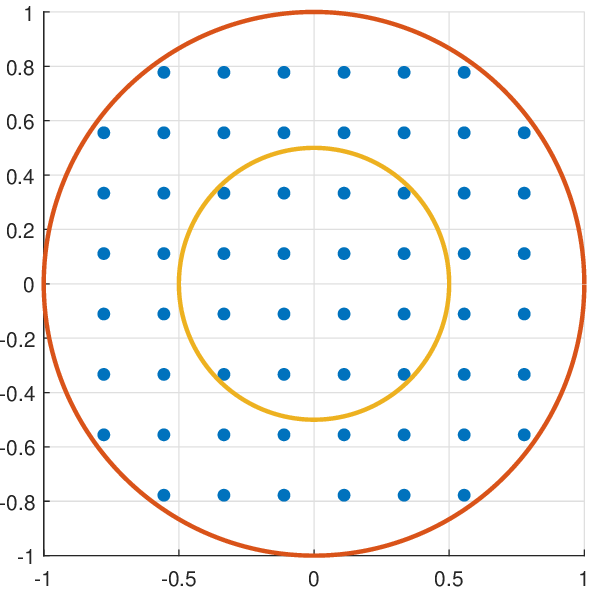}
	\caption{Example~\ref{exa2}: The sensors distribution. Blue dots represent the location of sensors, and solid lines are interface $\Gamma$ and $\partial \Omega$ in the case where $r_0=0.5$.}
	\label{Fig4}
\end{figure}

\begin{figure}[htbp]
	\centering
	\begin{minipage}{0.45\linewidth}
		\centering
		\includegraphics[width=1\linewidth]{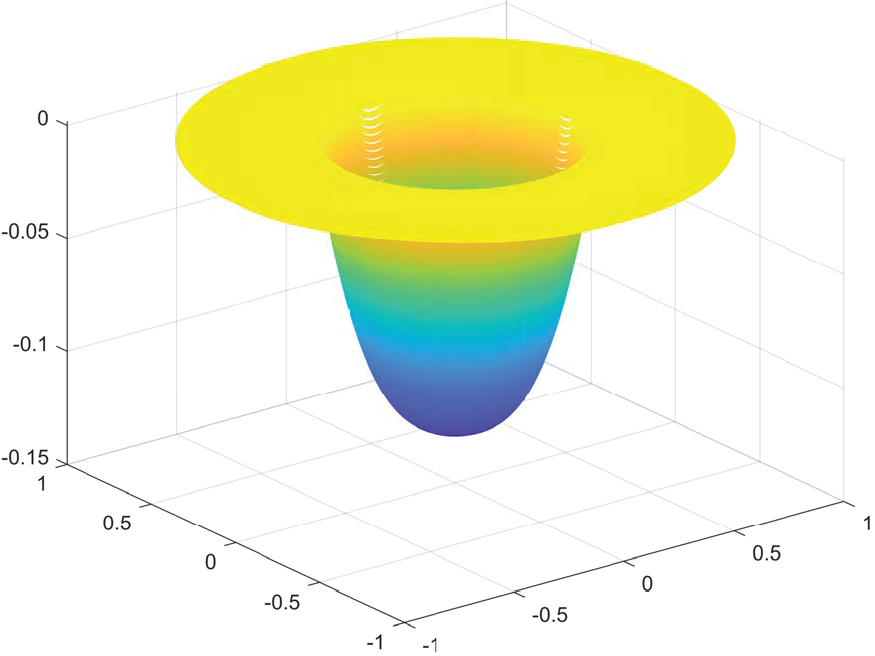}
	\end{minipage}
	\begin{minipage}{0.45\linewidth}
		\centering
		\includegraphics[width=1\linewidth]{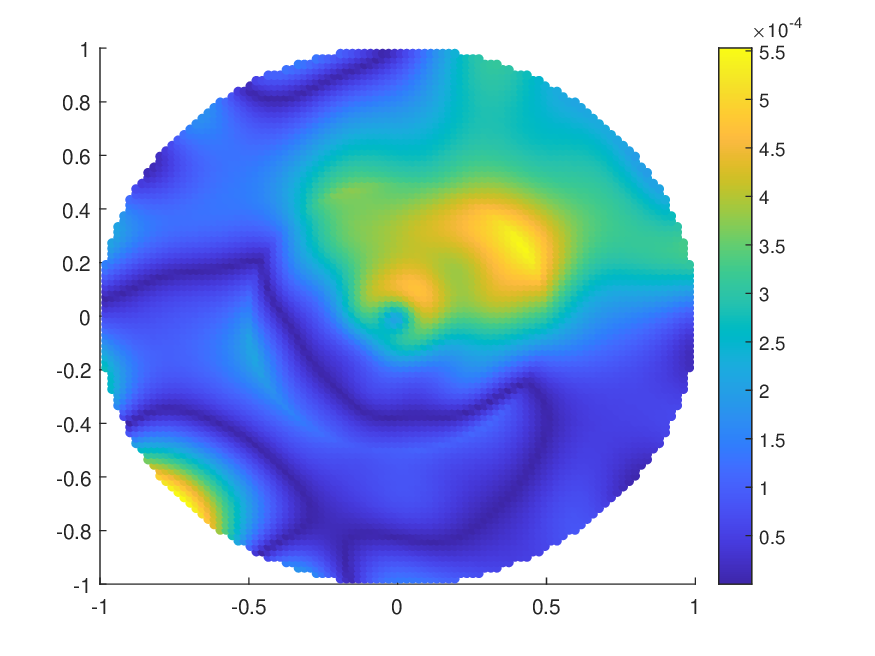}
	\end{minipage}
	\begin{minipage}{0.45\linewidth}
		\centering
		\includegraphics[width=1\linewidth]{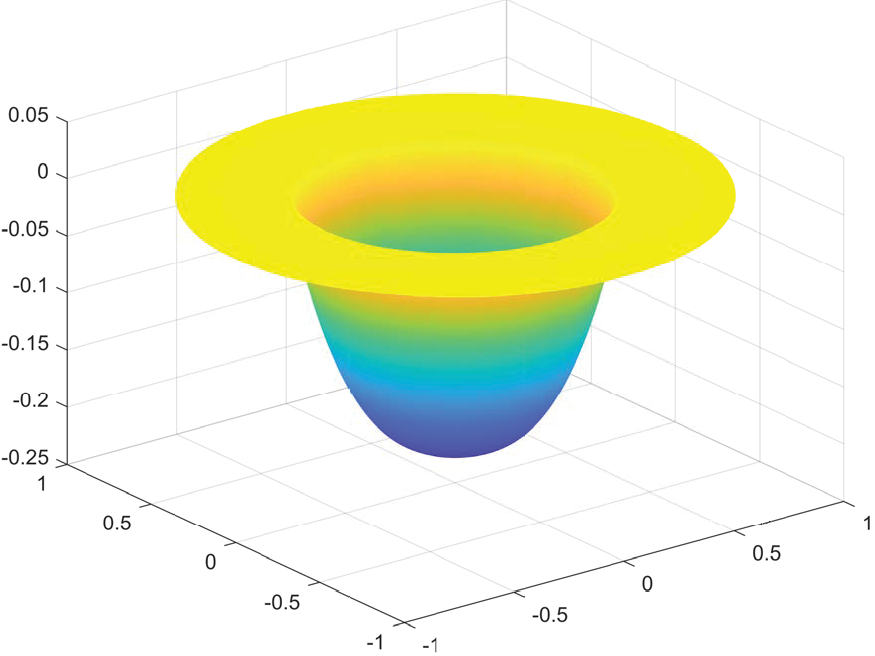}
	\end{minipage}
	\begin{minipage}{0.45\linewidth}
		\centering
		\includegraphics[width=1\linewidth]{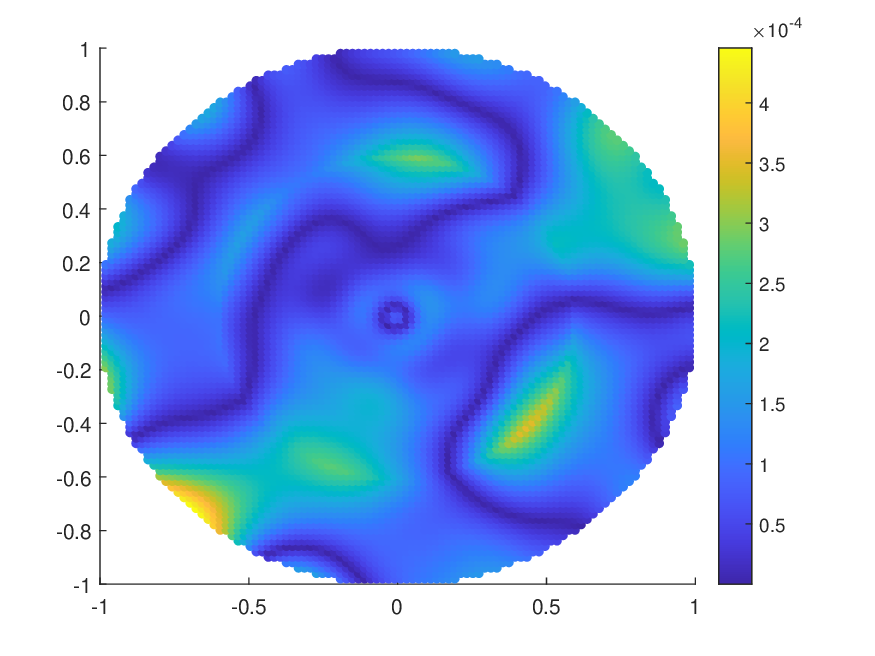}
	\end{minipage}
	\begin{minipage}{0.45\linewidth}
		\centering
		\includegraphics[width=1\linewidth]{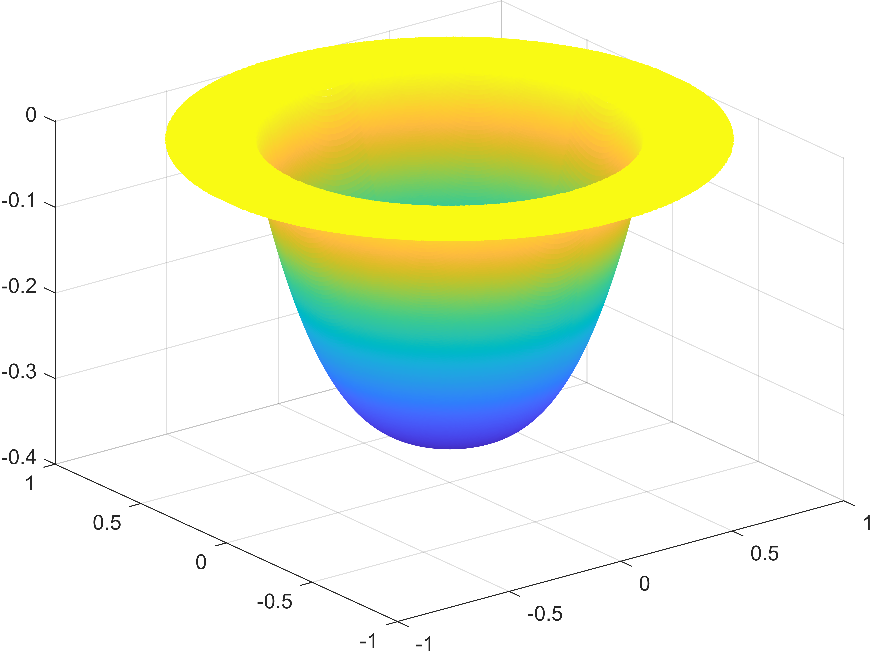}
	\end{minipage}
	\begin{minipage}{0.45\linewidth}
		\centering
		\includegraphics[width=1\linewidth]{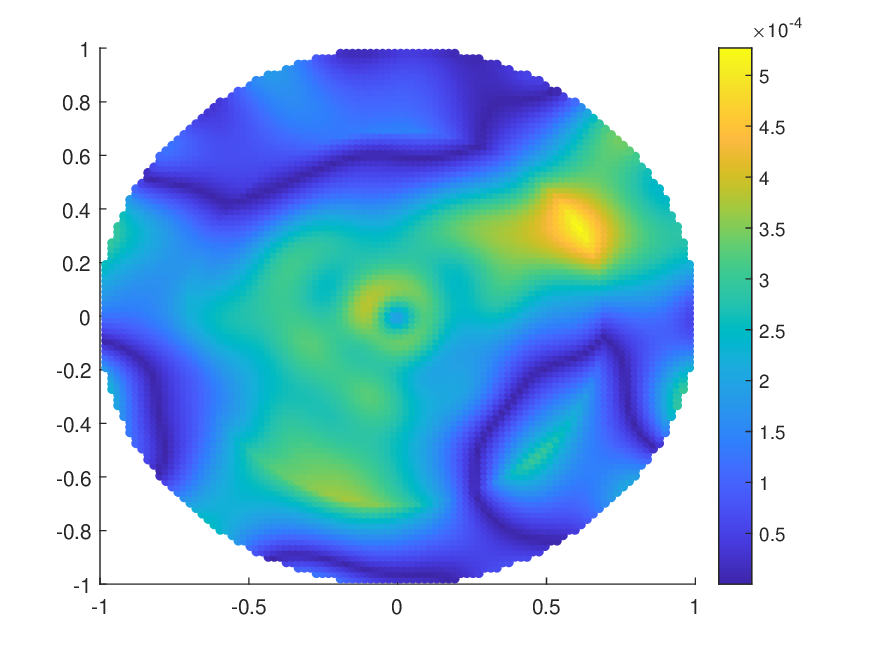}
	\end{minipage}
	\caption{Example~\ref{exa2}: The profile of NN solutions $u_{\mathcal{S}}$ and absolute point-wise errors  $\left|u - u_{\mathcal{S}}\right|$ for the considered problem with different interface position and size. Top: $r_0=0.5$. Middle: $r_0=0.6$. Bottom: $r_0=0.7$.}
	\label{Fig5}
\end{figure}

The positions of the 60 sensors in the input function are equidistant grid points in the computational domain $\Omega$ (see Fig.~\ref{Fig4} for an illustration). The number of input functions $(f, \Phi)$ is $500$. We separately calculate the relative $L^2$ errors between the exact and  approximation solutions for three cases of $r_0=0.5$, $r_0=0.6$ and $r_0=0.7$. As shown in Table~\ref{Table4} and Fig.~\ref{Fig5}, the XI-DeepONet performs very well for different interfaces. Especially, the relative $L^2$ error is about $1.30 \times 10^{-3}$ when $r_0=0.7$. These results confirm the capability of XI-DeepONet to solve high contrast diffusion coefficient interface problems with different size of the interface (here is the radius of interface circle).
\begin{table}[htbp]
	\centering
	\begin{tabular}{c|ccc}
		\toprule
		$r_0$    & No. of sensors  & No. of parameters & ${{{{\left\| {{u_S} - u} \right\|}_{L^2} }} \mathord{\left/{\vphantom {{{{\left\| {{u_S} - u} \right\|}_{L^2} }} {{{\left\| u \right\|}_{L^2} }}}} \right.\kern-\nulldelimiterspace} {{{\left\| u \right\|}_{L^2} }}}$\\
		\midrule
		$0.5$    &    60   &  133801 & $5.08 \times 10^{-3}$\\
		$0.6$    &    60   &  133801 & $1.48 \times 10^{-3}$\\
		$0.7$    &    60   &  133801 & $1.30 \times 10^{-3}$\\
		\bottomrule
	\end{tabular}%
	\caption{Example~\ref{exa2}: The relative $L^2$ errors of XI-DeepONet for different interface $\Gamma$($r_0=0.5, r_0=0.6, r_0=0.7$).}
	\label{Table4}
\end{table}

\begin{figure}[htbp]
	\centering
	\begin{minipage}{0.3\linewidth}
		\centering
		\includegraphics[width=1\linewidth]{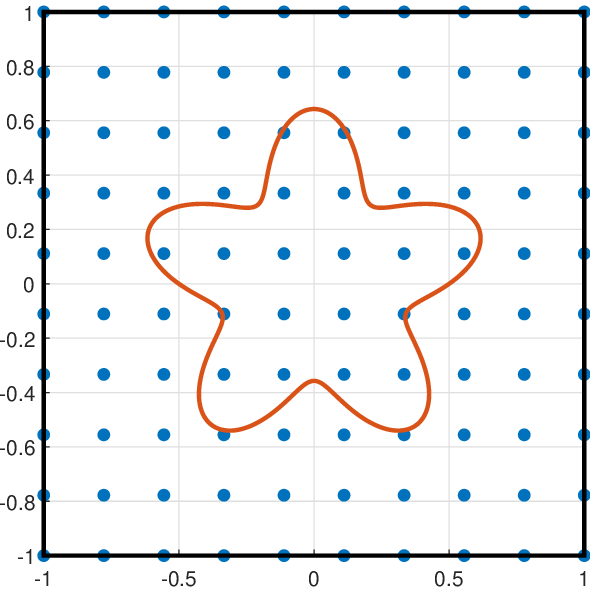}
	\end{minipage}
	\begin{minipage}{0.3\linewidth}
		\centering
		\includegraphics[width=1\linewidth]{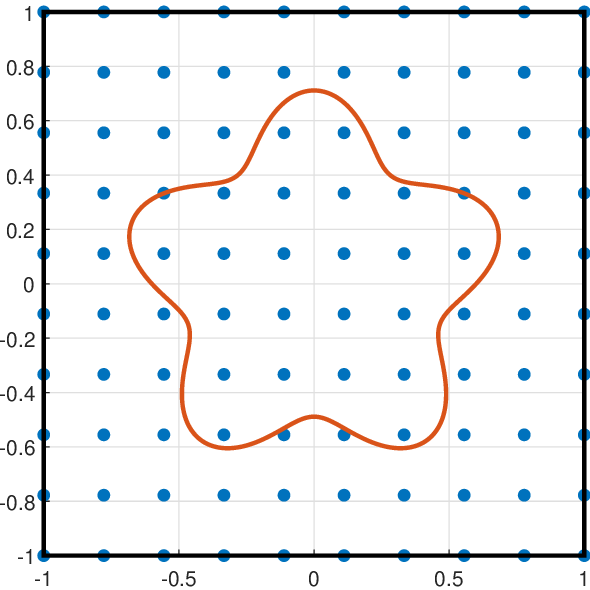}
	\end{minipage}
	\begin{minipage}{0.3\linewidth}
		\centering
		\includegraphics[width=1\linewidth]{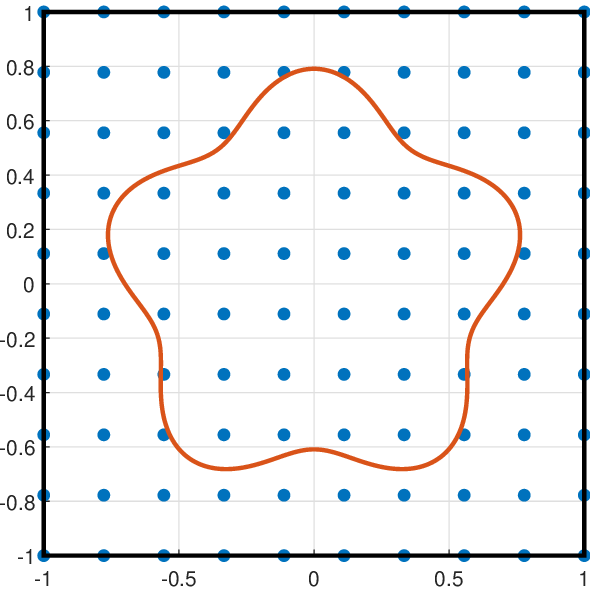}
	\end{minipage}
	\caption{Example~\ref{exa3}: The sensors distribution. Blue dots represent the location of sensors, and solid lines are interface $\Gamma$ and $\partial \Omega$.}
	\label{Fig6}
\end{figure}

\begin{example}\label{exa3}
	Consider the model problem \eqref{2.1}--\eqref{2.4} in the square $\Omega=[-1,1]^2$. 
	The diffusion coefficient is a piecewise constant, which is given by $a(x) = 2, x \in \Omega^+$ and  $a(x) = 1, x \in \Omega^-$. 
	The interface $\Gamma$ is irregular and complicated, which is given by the polar coordinate form as follows:
	\eq{\label{levelset2}\phi(r,\theta)=r-r_1-\dfrac{{\sin ( {5\theta } )}}{{{r_2}}},}
	where $(r_1, r_2)$ takes from $\left[0.5, 0.7\right] \times \left[ 7.0, 11.0\right]$ randomly. The interface conditions on $\Gamma$ are given as
	$$
	{g_D}( {x,y} ) = \dfrac{1}{{1 + 10( {{x^2}+{y^2}} )}}, \qquad
	{g_N}( {x,y} ) = 0,
	$$
	and the boundary condition is set as 
	$$
	h( {x,y} ) = \dfrac{2}{{1 + 10( {{x^2} + {y^2}} )}}.
	$$ 
\end{example}	
This example is taken from \cite{wu2024solving}, which is used to demonstrate the capability of XI-DeepONet to solve the parametric interface problem with irregular interfaces, which can be characterized by several parameters (here are $r_1$ and $r_2$ in \eqref{levelset2}). It is clear that we can obtain interfaces with different shapes by changing the value of $(r_1, r_2)$ (see Fig.~\ref{Fig6} for an illustration). Further, the source function $f$ takes the following formula:
$$
{f^\pm }( {x,y} ) = \dfrac{{p_1^\pm }}{{{{( {1 + 10( {{x^2} + {y^2}} )} )}^2}}} - \dfrac{{p_2^\pm ( {{x^2} + {y^2}} )}}{{{{( {1 + 10( {{x^2} + {y^2}} )} )}^3}}},
$$
where $( {p_1^ \pm },{p_2^ \pm })\in \left[ 50,100\right] \times \left[ 1550,1650\right]$. 
When $( {p_1^ + ,p_2^ + } ) = ( {p_1^ - ,p_2^ - } ) = ( {80,1600} )$, we can obtain the exact solution as follows:
\begin{align*}
u\left( x, y\right) =\left\{
\begin{aligned}
&\dfrac{2}{1+10\left( x^2 + y^2\right) } &\mathrm{in}\;\Omega^+,\\
&\dfrac{1}{1+10\left( x^2 + y^2\right) } &\mathrm{in}\;\Omega^-,
\end{aligned}
\right.
\end{align*}
which will be used as the exact solution for later testing. This example is an elliptic interface problem with non-homogeneous interface conditions. We can solve this problem using the XI-DeepONet framework, as mentioned in Remark~\ref{remark1}.

Unfortunately, the derivative of the level set function \eqref{levelset2} has a singularity at the origin. It is easy to see that 
$$
\mathop {\lim }\limits_{( {x,y} ) \to ( {0,0} )} \frac{{\partial \phi }}{{\partial x}} = \mathop {\lim }\limits_{( {x,y} ) \to ( {0,0} )} \frac{{\partial \phi }}{{\partial r}}\frac{{\partial r}}{{\partial x}} + \frac{{\partial \phi }}{{\partial \theta }}\frac{{\partial \theta }}{{\partial x}} = \infty.
$$
To solve this difficulty, we replace the absolute value function in \eqref{phifun} with the $ReLU$ function, i.e., 
$$\Phi (x,y) = ReLU( {\phi ( {x,y} )} ) = \left\{ {\begin{array}{*{20}{lc}}
	{\phi ( {x,y} ),}&{\mathrm{in} \; \Omega^+,}\\
	{0,}&{\mathrm{in} \; \Omega^-.}
	\end{array}} \right.$$
Then $\dfrac{{\partial \phi }}{{\partial x}} = \dfrac{{\partial \phi }}{{\partial y}} = 0\;\mathrm{in}\;\Omega^-$. 

\begin{figure}[htbp]
	\centering
	\begin{minipage}{0.45\linewidth}
		\centering
		\includegraphics[width=1\linewidth]{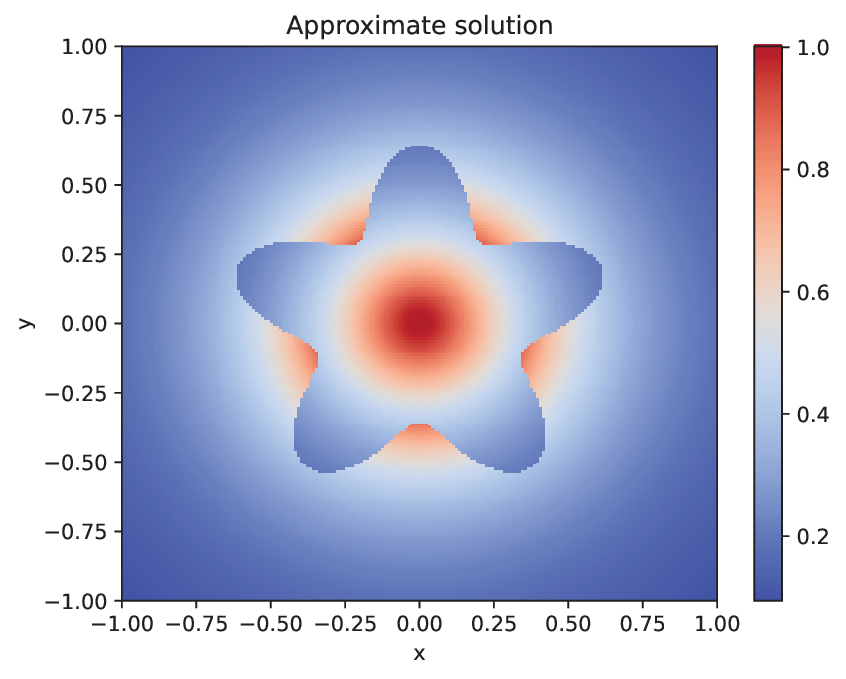}
	\end{minipage}
	\begin{minipage}{0.45\linewidth}
		\centering
		\includegraphics[width=1\linewidth]{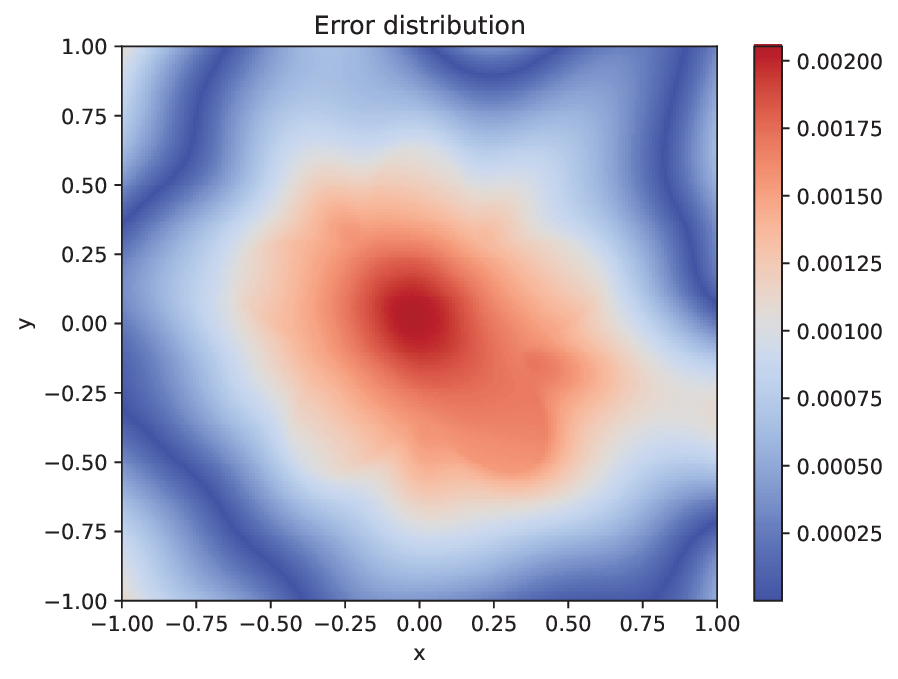}
	\end{minipage}
	\begin{minipage}{0.45\linewidth}
		\centering
		\includegraphics[width=1\linewidth]{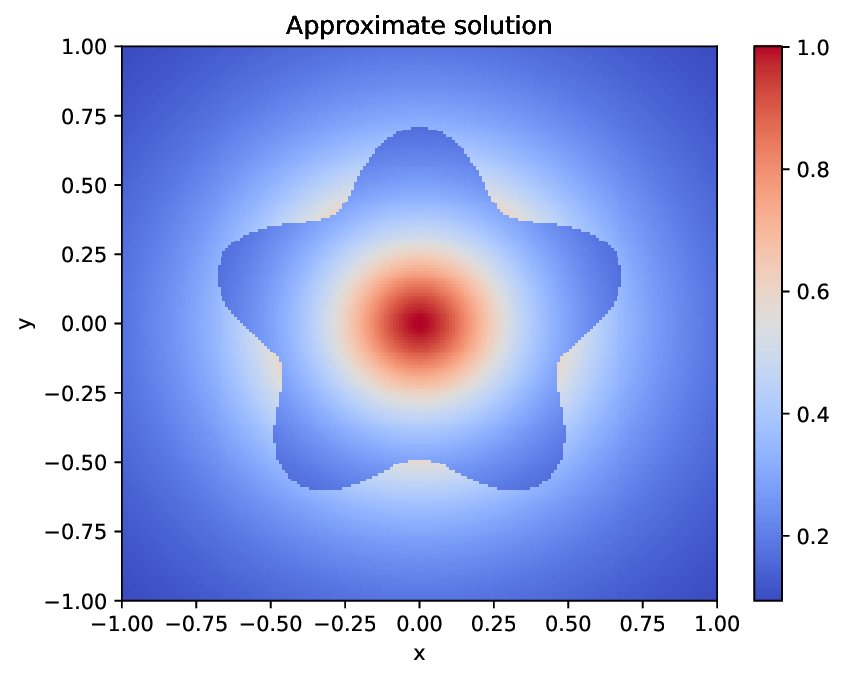}
	\end{minipage}
	\begin{minipage}{0.45\linewidth}
		\centering
		\includegraphics[width=1\linewidth]{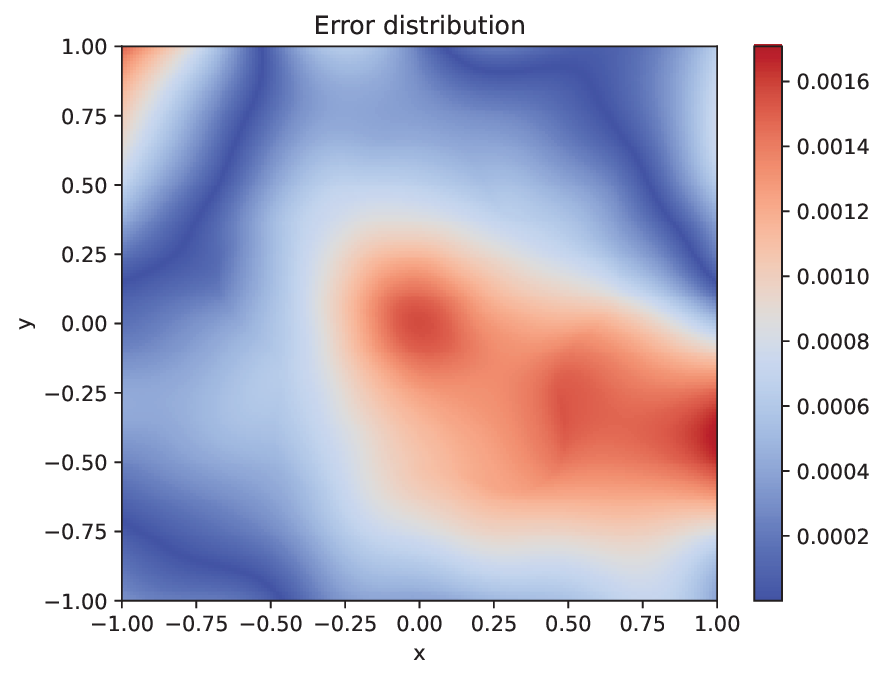}
	\end{minipage}
	\begin{minipage}{0.45\linewidth}
		\centering
		\includegraphics[width=1\linewidth]{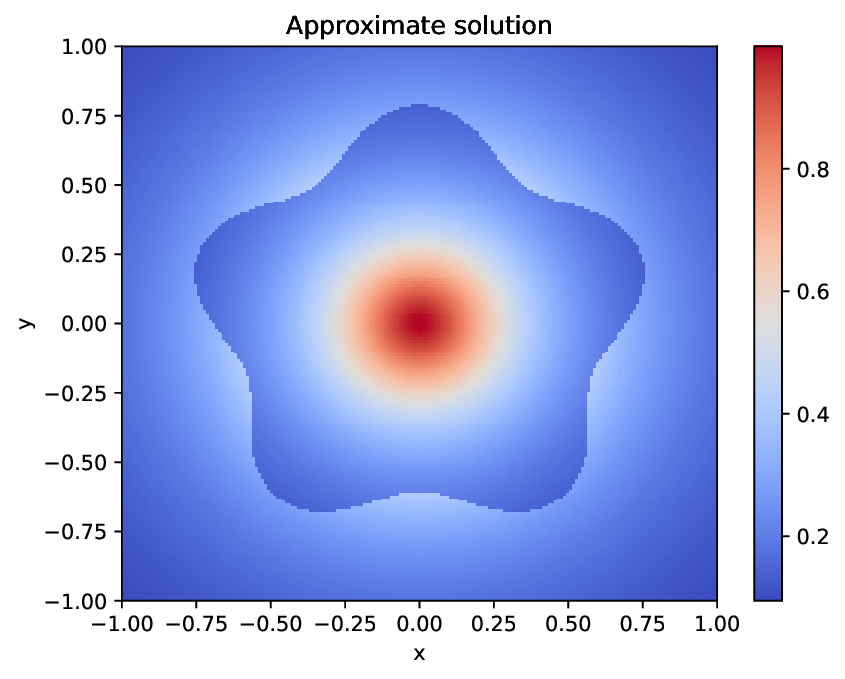}
	\end{minipage}
	\begin{minipage}{0.45\linewidth}
		\centering
		\includegraphics[width=1\linewidth]{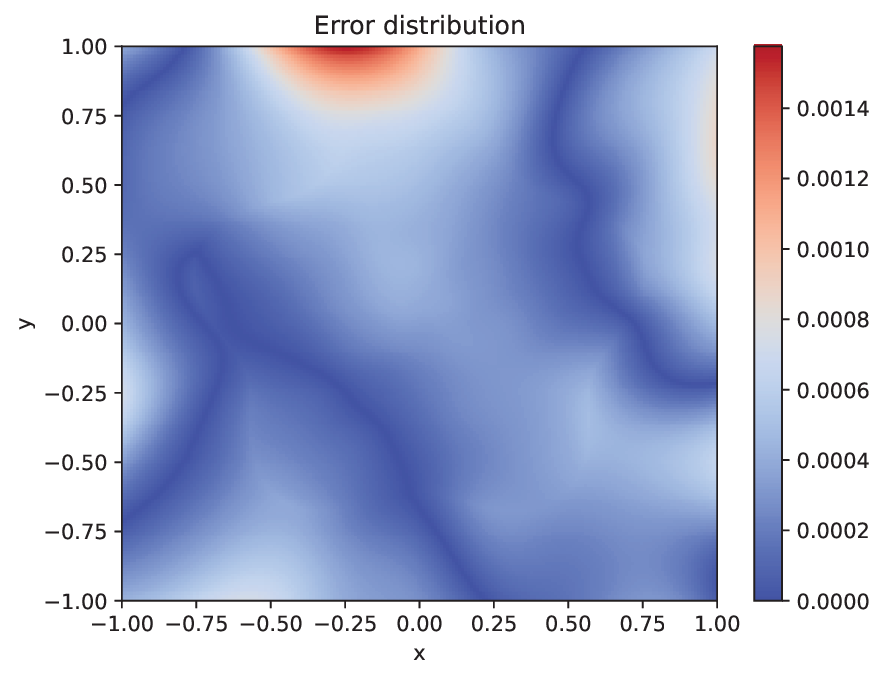}
	\end{minipage}
	\caption{Example~\ref{exa3}: The profile of network solutions $u_{\mathcal{S}}$ and absolute point-wise error  $\left|u - u_{\mathcal{S}}\right|$. Top: $(r_1, r_2 ) = (0.5, 7.0 )$. Middle: $(r_1, r_2 ) = (0.6, 9.0 )$. Bottom: $(r_1, r_2 ) = (0.7, 11.0 )$.}
	\label{Fig7}
\end{figure}

We randomly generate 500 input paired functions $(f,\Phi)$ to train the neural network. As shown in Fig.~\ref{Fig6}, there are 100 sensors uniformly distributed in the $\Omega$ for input functions. We set three level set function parameters, $(r_1, r_2 ) = (0.5, 7.0 )$, $(r_1, r_2 ) = (0.6, 9.0 )$ and $(r_1, r_2 ) = (0.7, 9.0 )$, to show the accuracy of the proposed NN method. As one can see from Table~\ref{Table5} and Fig.~\ref{Fig7}, the irregular interface problem can be handled properly by the XI-DeepONet model.
\begin{table}[htbp]
	\centering
	\begin{tabular}{c|ccc}
		\toprule
		$(r_1, r_2 ) $     & No. of sensors     &No. of parameters & ${{{{\left\| {{u_S} - u} \right\|}_{L^2} }} \mathord{\left/{\vphantom {{{{\left\| {{u_S} - u} \right\|}_{L^2} }} {{{\left\| u \right\|}_{L^2} }}}} \right.\kern-\nulldelimiterspace} {{{\left\| u \right\|}_{L^2} }}}$\\
		\midrule
		$(0.5, 7.0 ) $      & 100 &  141801   & $3.60 \times 10^{-3}$\\
		$(0.6, 9.0 ) $      & 100 &  141801   & $2.95 \times 10^{-3}$\\
		$(0.7, 11.0 ) $     & 100 &  141801   & $1.42 \times 10^{-3}$\\
		\bottomrule
	\end{tabular}%
	\caption{Example~\ref{exa3}: The relative $L^2$ errors for different interfaces.}
	\label{Table5}
\end{table}

\begin{example}\label{exa6}
	Consider the model problem \eqref{2.1}--\eqref{2.4} in the spherical shell 
	$$
	\Omega=\left\{ {\left. {( {x,y,z} ) \in {\mathbb{R}^3}} \right|{{0.151}^2} \le {x^2} + {y^2} + {z^2} \le {{0.911}^2}} \right\}.
	$$
	The diffusion coefficient is set to be
	$$
	a(x,y,z) = \left\{ 
	\begin{array}{lc}
	10 + 2\cos (2\pi (x + y) )\sin (2\pi(x - y))\cos ( z ) &{\mathrm{in} \; \Omega^-},\\
	1 &{\mathrm{in} \; \Omega^+}.
	\end{array}
	\right.
	$$
	The embedded interface $\Gamma$ is represented by the zero level set of function
	\[
	\begin{split}
	\phi (x,y,z) = & \sqrt{{x^2} + {y^2} + {z^2}} \\
	& -{r_0}\left( 1 + \left( \dfrac{{{x^2} + {y^2}}}{{{x^2} + {y^2} + {z^2}}} \right)^2\sum\limits_{k = 1}^3 t_k\cos\left( n_k( \arctan ( \dfrac{y}{x} ) - \theta _k )\right)  \right),
	\end{split}
	\]
where the parameters $r_0, t_k, n_k, \theta_k, k=1,2,3$ are real numbers randomly chosen from the following intervals:
	$$
	r_0 \in [0.45, 0.55],
	$$
	$$
	(t_1,t_2,t_3) \in [0,0.2]\times[-0.2,0]\times[0.1,0.2],
	$$
	$$
	(n_1,n_2,n_3) \in [2,4]\times[3,5]\times[6,8],
	$$
	$$
	(\theta_1,\theta_2,\theta_3) \in [0.3,0.7]\times[1.6,2.0]\times[-0.1,0.1].
	$$
	The illustration of the domain and interface geometry can be found in \cite{bochkov2020solving}.
	For later testing, the exact solution is taken to be
	\begin{equation*}
	u(x,y,z) = \left\{ 
	\begin{array}{lc}
	\sin (2x)\cos(2y){e^z} & \mathrm{in} \; \Omega^-,\\
	16\Big(\dfrac{y - x}{3}\Big)^5 - 20\Big(\dfrac{y - x}{3}\Big)^3 + 5\Big(\dfrac{y - x}{3}\Big) & \\
	\qquad + \log ( {x + y + 3} )\cos ( z )& {\mathrm{in}\; \Omega^+},
	\end{array}
	\right.
	\end{equation*}
	and the corresponding source term $f$ takes the following formula:
	\begin{equation}
	f(x,y,z) = - \nabla  \cdot ( {a(x, y, z) \nabla u(x, y, z)} ) + \epsilon. \nonumber
	\end{equation}
	where $\epsilon \in \left[-0.1, 0.1\right]$ is a random number. When $\epsilon=0$ we can obtain the exact solution $u$. 
\end{example}

Solving three-dimensional elliptic interface problems with irregular interfaces $\Gamma$ poses significant challenges for traditional numerical methods. The cost is further compounded when the task involves solving hundreds or even thousands of such problems. Additionally, the computational expense associated with training a multitude of models for varying interface problems renders PINNs less feasible. This example aims to investigate the efficacy of the XI-DeepONet approach in tackling 3-D complex interface problems ( cf. \cite{bochkov2020solving, hu2022discontinuity}). 

\begin{figure}[htbp]
	\label{Fig9}
	\centering
	\begin{minipage}{0.45\linewidth}
		\centering
		\includegraphics[width=1\linewidth]{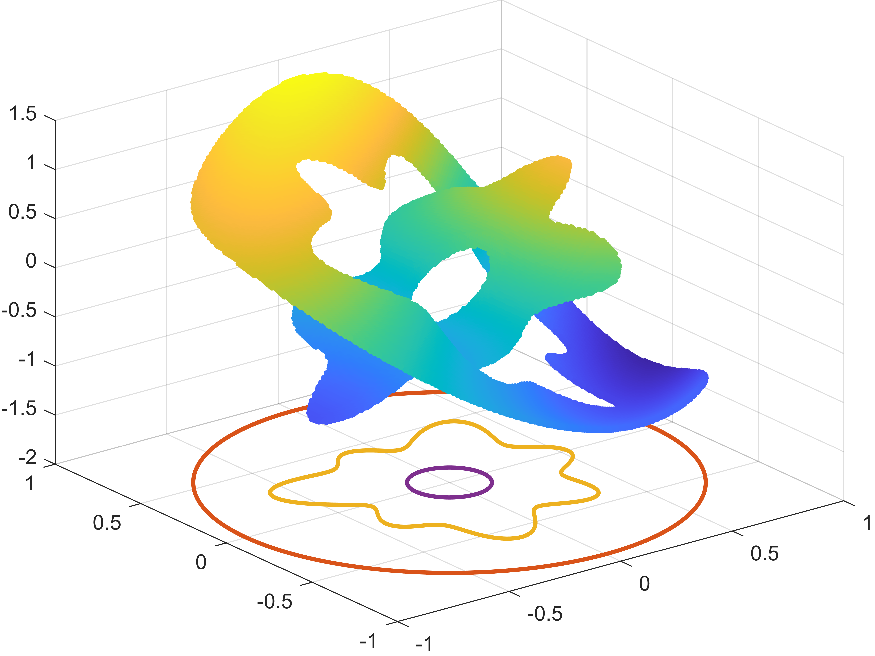}
	\end{minipage}
	\begin{minipage}{0.45\linewidth}
		\centering
		\includegraphics[width=1\linewidth]{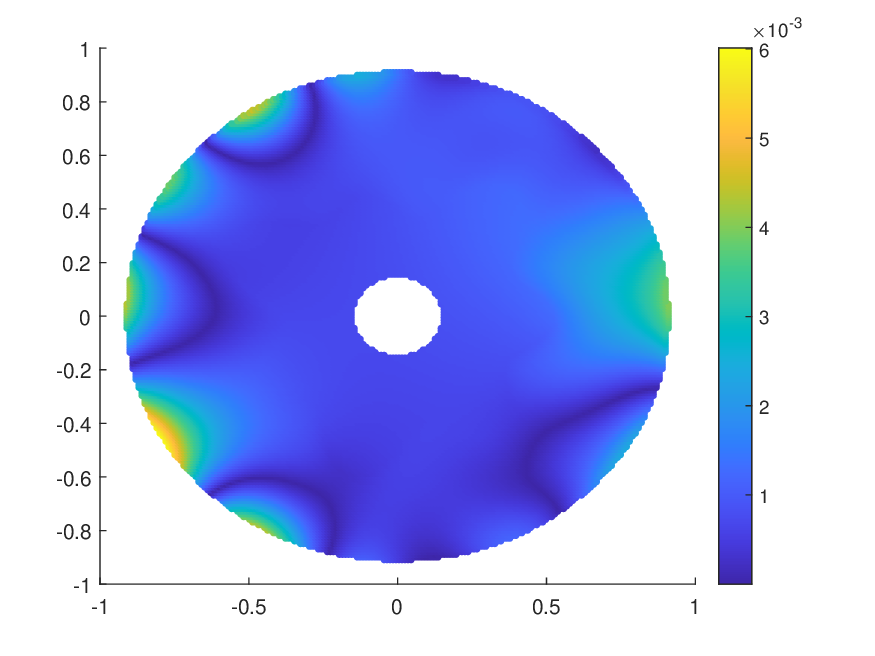}
	\end{minipage}
	\begin{minipage}{0.45\linewidth}
		\centering
		\includegraphics[width=1\linewidth]{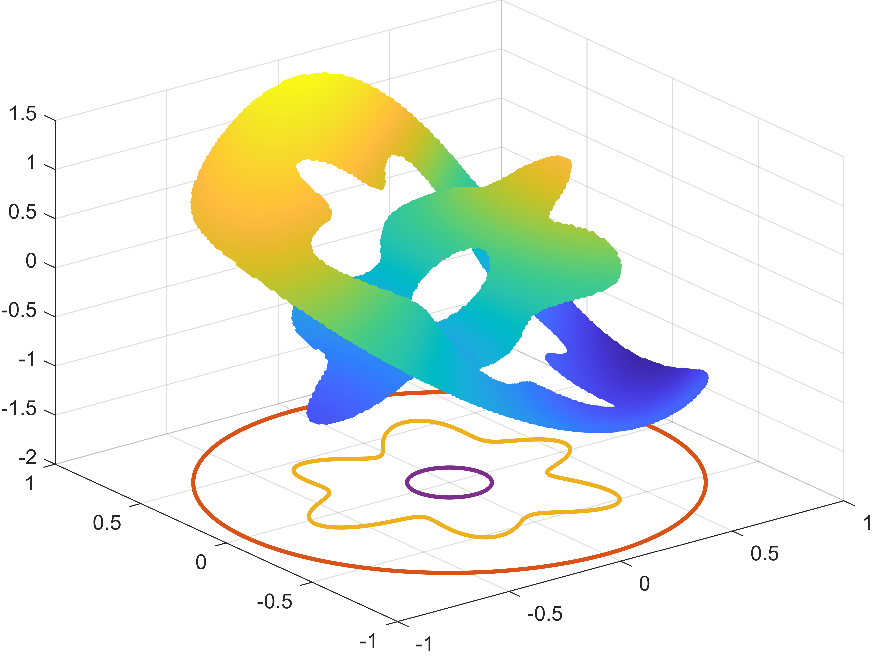}
	\end{minipage}
	\begin{minipage}{0.45\linewidth}
		\centering
		\includegraphics[width=1\linewidth]{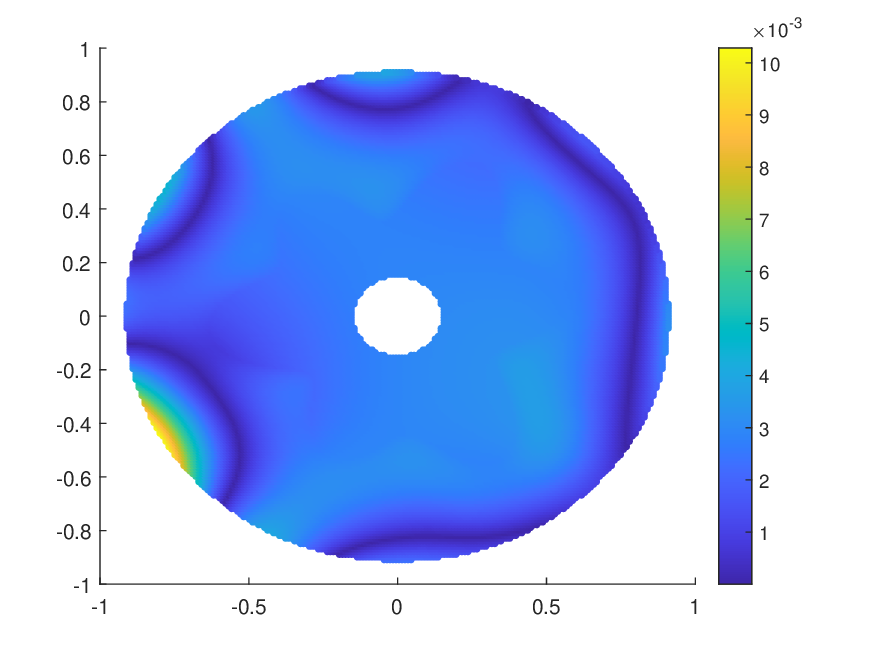}
	\end{minipage}
	\begin{minipage}{0.45\linewidth}
		\centering
		\includegraphics[width=1\linewidth]{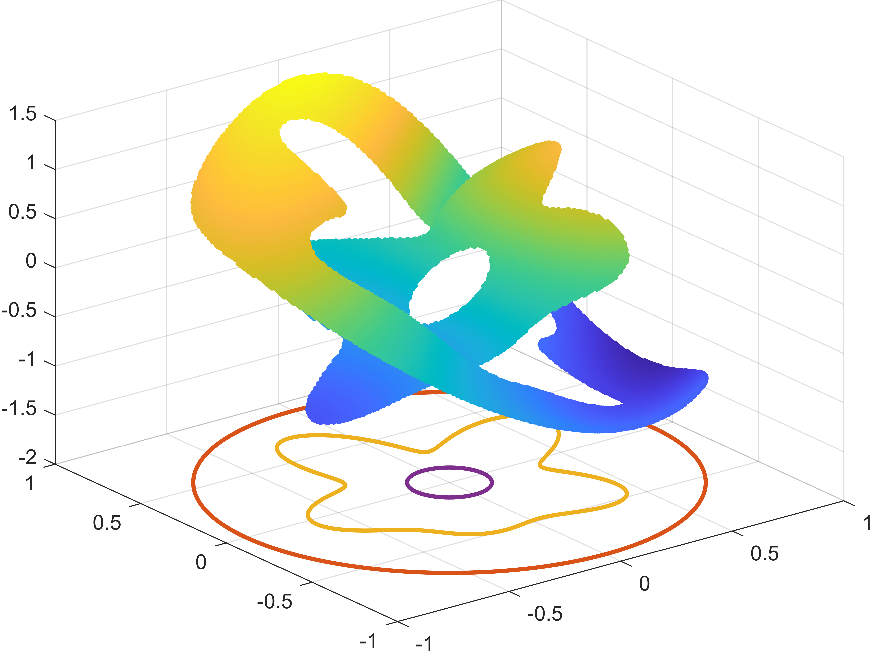}
	\end{minipage}
	\begin{minipage}{0.45\linewidth}
		\centering
		\includegraphics[width=1\linewidth]{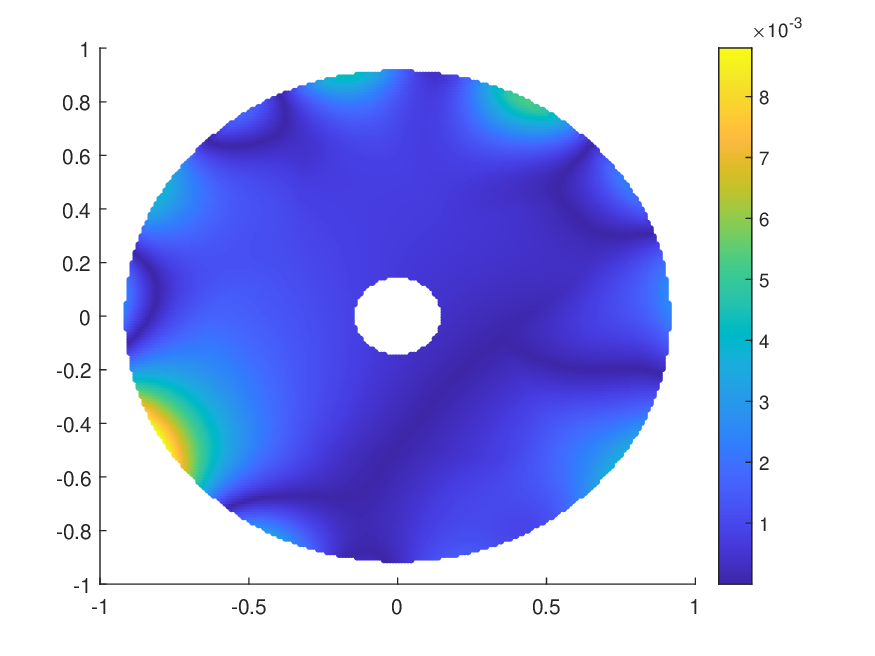}
	\end{minipage}
	\caption{The cross-sectional profile of $u_\mathcal{S}$ and absolute error on $z=0$ in Example~\ref{exa6}. Top: $\Gamma_1$, Middle: $\Gamma_2$, Bottom: $\Gamma_3$. The curves represent corresponding cross-sectional interface and boundary.}
\end{figure}

In this example, we learn the solution operator $\mathcal{G}$ mapping from $\Phi$ and $f$ to the latent solution $u$. We sampled 500 input level set functions $\phi(x,y,z)$ for training networks, and choose three sets of parameters to test the accuracy of results. Here, the branch and truck networks are 5-layer FNNs with 150 units in each hidden layer and there are 136 sensors distributed in computational domain $\Omega$. The setup of three testing interface $\Gamma_i,i=1,2,3$ and the relative $L^2$ errors with $\Gamma_i$ are shown in Table \ref{Table6}. Fig. \ref{Fig9} shows the cross-sectional profile of the network solution and absolute point-wise error on the hyperplane $z=0$. It is shown that the XI-DeepONets can be used to solve 3-D elliptic interface problems in the domain with complex interface problems subject to nonzero solution jump conditions. More importantly, we only need to train the neural networks once.
\begin{table}[htbp]
	\label{Table6}
	\centering
	\begin{tabular}{c|ccccc}
		\toprule
		& $r_0$    & $( t_1,t_2,t_3) $  &$( n_1,n_2,n_3)$ &$ ( \theta_1, \theta_2, \theta_3 ) $&  $\dfrac{\|u_S - u\|_{L^2}}{\|u\|_{L^2} }$\\
		\midrule
		$\Gamma_1$     &    0.483   & (0.10,-0.10,0.15) &  (3.00, 4.00, 7.00)  & (0.50, 1.80, 0.00) & $2.23 \times 10^{-3}$\\
		$\Gamma_2$      &    0.500   & (0.00,-0.10,0.20) &  (3.00, 3.00, 6.00)  &  (0.60, 1.90, 0.10) & $4.60 \times 10^{-3}$\\
		$\Gamma_3$      &    0.530   & (0.20, 0.20, 0.00) & (4.00, 5.00, 8.00) &  (0.40, 1.80, 0.10) & $2.28 \times 10^{-3}$\\
		\bottomrule
	\end{tabular}%
	\caption{Example~\ref{exa6}: The relative $L^2$ errors for different interface parameters.}
\end{table}

As the last example, we want to investigate the ability to solve high-dimensional elliptic interface problems by taking the dimension $d=6$. The problem is setup as follows \cite{lai2022shallow}. 
\begin{example}\label{exa5}
	Consider a 6-sphere of radius $0.6$ as the domain $\Omega$ enclosing another smaller 6-sphere of radius $r_0$ as the interior region $\Omega^-$, where $r_0 \in [0.4,0.5]$. The spherical interface can be labeled a zero-level set of the function
	$$\phi ( x ) = \left\| x \right\|_2^2 - {r_0}^2.$$
	The exact solution \cite{lai2022shallow} is chosen as
	\begin{equation}
	u( {x} ) = \left\{ {\begin{array}{*{20}{lc}}
		{\exp ( {{r_0^2} - \left\| {x} \right\|_2^2} ) + \sum\limits_{k = 1}^5 {\sin ( {{x_k}} )} }&{\mathrm{in} \; \Omega^+},\\
		{1 + 2\sin ( {{r_0^2} - \left\| {x} \right\|_2^2} ) + \sum\limits_{k = 1}^5 {\sin ( {{x_k}} )} }&{\mathrm{in} \; \Omega^-}.
		\end{array}} \right.
	\nonumber
	\end{equation}
\end{example}

The source term and interface conditions can be obtained using equation (\ref{2.1})--(\ref{2.4}). We generate 233 uniformly distributed sensors within domain $\Omega$. The results are shown in Table \ref{Table7}, the relative $L^2$ errors are about $7 \times 10^{-4}$. This example shows that the present method is applicable to high-dimensional elliptic interface problems.
\begin{table}[htbp]
	\centering
	\caption{Example~\ref{exa5}: The relative $L^2$ errors with different spherical interface radius.}
	\label{Table7}
	\begin{tabular}{c|cccc}
		\toprule
		$r_0$     & No. of sensors  & No. of parameters & ${{{{\left\| {{u_S} - u} \right\|}_{L^2} }} \mathord{\left/{\vphantom {{{{\left\| {{u_S} - u} \right\|}_{L^2} }} {{{\left\| u \right\|}_{L^2} }}}} \right.\kern-\nulldelimiterspace} {{{\left\| u \right\|}_{L^2} }}}$\\
		\midrule
		$0.4$    &  233   & 168801   & $7.62 \times 10^{-4}$\\
		$0.45$   &  233   & 168801   & $6.18 \times 10^{-4}$\\
		$0.5$    &  233   & 168801   & $7.14 \times 10^{-4}$\\
		\bottomrule
	\end{tabular}%
\end{table}

\section{Conclusion}
\label{Section 5}
We have proposed an extended interface deep operator network to solve parametric elliptic interface problems. This new neural operator can take the interface represented by the level set function as the input feature. Therefore, it has the generalization ability to interfaces which may have different positions and shapes. Furthermore, the PI-DeepONet framework is adopted to effectively reduce the demand for training datasets. Numerous numerical examples including complicated interface and high-dimension problem are carried out to show the efficiency of the proposed method. It is found that the new NN method can capture the discontinuities of solutions across the interface very efficiently. 


\bibliography{arxiv}

\begin{thebibliography}{10}
\expandafter\ifx\csname url\endcsname\relax
  \def\url#1{\texttt{#1}}\fi
\expandafter\ifx\csname urlprefix\endcsname\relax\def\urlprefix{URL }\fi
\expandafter\ifx\csname href\endcsname\relax
  \def\href#1#2{#2} \def\path#1{#1}\fi

\bibitem{greengard1994numerical}
L.~Greengard, M.~Moura, On the numerical evaluation of electrostatic fields in composite materials, Acta numerica 3 (1994) 379--410.

\bibitem{stroud2002numerical}
J.~Stroud, S.~Berger, D.~Saloner, Numerical analysis of flow through a severely stenotic carotid artery bifurcation, Journal of Biomechanical Engineering 124~(1) (2002) 9--20.

\bibitem{zienkiewicz2005finite}
O.~C. Zienkiewicz, R.~L. Taylor, J.~Z. Zhu, The finite element method: its basis and fundamentals, Elsevier, 2005.

\bibitem{mu2013weak}
L.~Mu, J.~Wang, G.~Wei, X.~Ye, S.~Zhao, Weak galerkin methods for second order elliptic interface problems, Journal of computational physics 250 (2013) 106--125.

\bibitem{chen1998finite}
Z.~Chen, J.~Zou, Finite element methods and their convergence for elliptic and parabolic interface problems, Numerische Mathematik 79~(2) (1998) 175--202.

\bibitem{li2010optimal}
J.~Li, J.~M. Melenk, B.~Wohlmuth, J.~Zou, Optimal a priori estimates for higher order finite elements for elliptic interface problems, Applied numerical mathematics 60~(1-2) (2010) 19--37.

\bibitem{guyomarc2009discontinuous}
G.~Guyomarc'h, C.-O. Lee, K.~Jeon, A discontinuous galerkin method for elliptic interface problems with application to electroporation, Communications in numerical methods in engineering 25~(10) (2009) 991--1008.

\bibitem{burman2010interior}
E.~Burman, P.~Hansbo, Interior-penalty-stabilized lagrange multiplier methods for the finite-element solution of elliptic interface problems, IMA journal of numerical analysis 30~(3) (2010) 870--885.

\bibitem{dryja2003discontinuous}
M.~Dryja, On discontinuous galerkin methods for elliptic problems with discontinuous coefficients, Computational Methods in Applied Mathematics 3~(1) (2003) 76--85.

\bibitem{he2020interface}
X.~He, W.~Deng, H.~Wu, An interface penalty finite element method for elliptic interface problems on piecewise meshes, Journal of Computational and Applied Mathematics 367 (2020) 112473.

\bibitem{leveque1994immersed}
R.~J. LeVeque, Z.~Li, The immersed interface method for elliptic equations with discontinuous coefficients and singular sources, SIAM Journal on Numerical Analysis 31~(4) (1994) 1019--1044.

\bibitem{li1998immersed}
Z.~Li, The immersed interface method using a finite element formulation, Applied Numerical Mathematics 27~(3) (1998) 253--267.

\bibitem{li2003new}
Z.~Li, T.~Lin, X.~Wu, New cartesian grid methods for interface problems using the finite element formulation, Numerische Mathematik 96 (2003) 61--98.

\bibitem{gong2008immersed}
Y.~Gong, B.~Li, Z.~Li, Immersed-interface finite-element methods for elliptic interface problems with nonhomogeneous jump conditions, SIAM Journal on Numerical Analysis 46~(1) (2008) 472--495.

\bibitem{gong2010immersed}
Y.~Gong, Z.~Li, Immersed interface finite element methods for elasticity interface problems with non-homogeneous jump conditions, Numerical Mathematics: Theory, Methods and Applications 3~(1) (2010) 23--39.

\bibitem{lin2019nonconforming}
T.~Lin, D.~Sheen, X.~Zhang, A nonconforming immersed finite element method for elliptic interface problems, Journal of Scientific Computing 79 (2019) 442--463.

\bibitem{hansbo2002unfitted}
A.~Hansbo, P.~Hansbo, An unfitted finite element method, based on nitsche’s method, for elliptic interface problems, Computer methods in applied mechanics and engineering 191~(47-48) (2002) 5537--5552.

\bibitem{xiao2020high}
Y.~Xiao, J.~Xu, F.~Wang, High-order extended finite element methods for solving interface problems, Computer Methods in Applied Mechanics and Engineering 364 (2020) 112964.

\bibitem{zhang2022condensed}
Q.~Zhang, C.~Cui, U.~Banerjee, I.~Babu{\v{s}}ka, A condensed generalized finite element method (cgfem) for interface problems, Computer Methods in Applied Mechanics and Engineering 391 (2022) 114537.

\bibitem{raissi2019physics}
M.~Raissi, P.~Perdikaris, G.~E. Karniadakis, Physics-informed neural networks: A deep learning framework for solving forward and inverse problems involving nonlinear partial differential equations, Journal of Computational physics 378 (2019) 686--707.

\bibitem{yu2018deep}
W.~E, B.~Yu, The deep ritz method: a deep learning-based numerical algorithm for solving variational problems, Communications in Mathematics and Statistics 6~(1) (2018) 1--12.

\bibitem{chen2022rfm}
J.~Chen, X.~Chi, W.~E, Z.~Yang, Bridging traditional and machine learning-based algorithms for solving pdes: The random feature method, Journal of Machine Learning 1~(3) (2022) 268--298.

\bibitem{lai2022shallow}
M.-C. Lai, C.-C. Chang, W.-S. Lin, W.-F. Hu, T.-S. Lin, A shallow ritz method for elliptic problems with singular sources, Journal of Computational Physics 469 (2022) 111547.

\bibitem{wu2022inn}
S.~Wu, B.~Lu, Inn: Interfaced neural networks as an accessible meshless approach for solving interface pde problems, Journal of Computational Physics 470 (2022) 111588.

\bibitem{he2022mesh}
C.~He, X.~Hu, L.~Mu, A mesh-free method using piecewise deep neural network for elliptic interface problems, Journal of Computational and Applied Mathematics 412 (2022) 114358.

\bibitem{hu2022discontinuity}
W.-F. Hu, T.-S. Lin, M.-C. Lai, A discontinuity capturing shallow neural network for elliptic interface problems, Journal of Computational Physics 469 (2022) 111576.

\bibitem{tseng2023cusp}
Y.-H. Tseng, T.-S. Lin, W.-F. Hu, M.-C. Lai, A cusp-capturing pinn for elliptic interface problems, Journal of Computational Physics 491 (2023) 112359.

\bibitem{chi2024random}
X.~Chi, J.~Chen, Z.~Yang, The random feature method for solving interface problems, Computer Methods in Applied Mechanics and Engineering 420 (2024) 116719.

\bibitem{lu2021learning}
L.~Lu, P.~Jin, G.~Pang, Z.~Zhang, G.~E. Karniadakis, Learning nonlinear operators via deeponet based on the universal approximation theorem of operators, Nature machine intelligence 3~(3) (2021) 218--229.

\bibitem{li2020fourier}
Z.~Li, N.~Kovachki, K.~Azizzadenesheli, B.~Liu, K.~Bhattacharya, A.~Stuart, A.~Anandkumar, Fourier neural operator for parametric partial differential equations, arXiv preprint arXiv:2010.08895 (2020).

\bibitem{lu2022comprehensive}
L.~Lu, X.~Meng, S.~Cai, Z.~Mao, S.~Goswami, Z.~Zhang, G.~E. Karniadakis, A comprehensive and fair comparison of two neural operators (with practical extensions) based on fair data, Computer Methods in Applied Mechanics and Engineering 393 (2022) 114778.

\bibitem{kontolati2024learning}
K.~Kontolati, S.~Goswami, G.~Em~Karniadakis, M.~D. Shields, Learning nonlinear operators in latent spaces for real-time predictions of complex dynamics in physical systems, Nature Communications 15~(1) (2024) 5101.

\bibitem{tripura2023wavelet}
T.~Tripura, S.~Chakraborty, Wavelet neural operator for solving parametric partial differential equations in computational mechanics problems, Computer Methods in Applied Mechanics and Engineering 404 (2023) 115783.

\bibitem{wang2021learning}
S.~Wang, H.~Wang, P.~Perdikaris, Learning the solution operator of parametric partial differential equations with physics-informed deeponets, Science advances 7~(40) (2021) eabi8605.

\bibitem{li2024physics}
Z.~Li, H.~Zheng, N.~Kovachki, D.~Jin, H.~Chen, B.~Liu, K.~Azizzadenesheli, A.~Anandkumar, Physics-informed neural operator for learning partial differential equations, ACM/JMS Journal of Data Science 1~(3) (2024) 1--27.

\bibitem{jin2022mionet}
P.~Jin, S.~Meng, L.~Lu, Mionet: Learning multiple-input operators via tensor product, SIAM Journal on Scientific Computing 44~(6) (2022) A3490--A3514.

\bibitem{wu2024solving}
S.~Wu, A.~Zhu, Y.~Tang, B.~Lu, Solving parametric elliptic interface problems via interfaced operator network, Journal of Computational Physics (2024) 113217.

\bibitem{baydin2018automatic}
A.~G. Baydin, B.~A. Pearlmutter, A.~A. Radul, J.~M. Siskind, Automatic differentiation in machine learning: a survey, Journal of machine learning research 18~(153) (2018) 1--43.

\bibitem{kingma2014adam}
D.~P. Kingma, J.~Ba, Adam: A method for stochastic optimization, arXiv preprint arXiv:1412.6980 (2014).

\bibitem{yu2007matched}
S.~Yu, Y.~Zhou, G.-W. Wei, Matched interface and boundary (mib) method for elliptic problems with sharp-edged interfaces, Journal of Computational Physics 224~(2) (2007) 729--756.

\bibitem{bochkov2020solving}
D.~Bochkov, F.~Gibou, Solving elliptic interface problems with jump conditions on cartesian grids, Journal of Computational Physics 407 (2020) 109269.

\end{thebibliography}
\end{document}